\documentclass[11pt]{article}
\usepackage{amsfonts,amssymb,amsmath,pdfsync,color,bbm,ulem,cancel}
\hsize \textwidth
\advance \hsize by -\marginparwidth
\evensidemargin \oddsidemargin
\usepackage{amssymb}
\advance\hoffset by 5mm

\newcommand{\Rm}{{\mathbb R}}
\newcommand{\Sm}{{\mathbb S}}
\newcommand{\eps}{\varepsilon}
\newcommand{\commentout}[1]{}

\newcommand{\bal}{\begin{aligned}}
\newcommand{\enbal}{\end{aligned}}
\newcommand{\be}{\begin{equation}}
\newcommand{\ee}{\end{equation}}

\newcommand{\disp}{\displaystyle}

\newcommand{\cT}{\mathcal{T}}
\newcommand{\Em}{{\mathbb E}}

\DeclareMathOperator{\pbmo}{p-BMO}

\def\red{\textcolor{red}}

\usepackage[colorlinks=true,pdfpagemode=UseNone,urlcolor=blue,linkcolor=blue,citecolor=blue,breaklinks=true]{hyperref} 

\newcommand{\one}{\mathbbm{1}}

\newtheorem{thm}{Theorem}[section]
\newtheorem{lem}[thm]{Lemma}
\newtheorem{cor}[thm]{Corollary}
\newtheorem{prop}[thm]{Proposition}

 \def\di{\displaystyle}

\newcommand{\farc}{\frac}
\numberwithin{equation}{section}
\numberwithin{thm}{section}

\begin{document}

\title{A uniform bound for solutions to a thermo-diffusive system}

\author{Joonhyun La \and Jean-Michel Roquejoffre \and Lenya Ryzhik}

\maketitle

\begin{abstract}
We obtain uniform in time $L^\infty$-bounds for the solutions to a class of thermo-diffusive systems.
The nonlinearity is assumed to be at most sub-exponentially growing at infinity and have a linear behavior near zero. 
\end{abstract}

\section{Introduction}

\subsubsection*{The thermo-diffusive system}

We investigate the global boundedness of the positive solutions to the initial
value problem for reaction-diffusion systems of the form
\begin{equation}\label{oct1702}
\begin{aligned}
&\partial_t u= \nu\Delta u - ug (v),\\
&\partial_t v = \kappa \Delta v + ug (v),
\end{aligned}  
\end{equation}
posed in the whole space $x\in\Rm^n$, with the nonlinear term $g(v)$ 
having possibly super-linear growth. It  has been known since Fujita's work \cite{F1966}
that solutions to the nonlinear heat equation
\[
\partial_tv=\Delta v+v^p
\]
may blow up in finite time. On the other hand, in the special case $\kappa=\nu$ of the
system (\ref{oct1702}), the function $z(t,x)=u(t,x)+v(t,x)$ solves the linear heat equation
\[
\partial_tz=\kappa\Delta z,
\]
so that both $u(t,x)$ and $v(t,x)$ are globally bounded in time. In this work,
we address the question of how the interplay between the two equations
in (\ref{oct1702}) may lead to uniform boundedness
of the solutions. 

In (\ref{oct1702}), we think of $u(t,x)$ as fuel concentration that is consumed
in the chemical reaction, and of~$v(t,x)$ as the concentration of the products of the reaction, or as the reactant temperature.
The systems we discuss in this paper
appear naturally, for example, in the exothermic combustion modeling and ecology \cite{CS1977,MS1974,MR2005,OALD1997}. 
  In accordance with 
this interpretation, we will assume that the reaction rate $g(v)$ is positive:
\be\label{feb102}
g(v)\ge 0,~~\hbox{ for all $v\ge 0$.}
\ee
This assumption, together with the comparison principle, ensures that 
\be\label{feb104}
\hbox{$u(t,x)> 0$ and $v(t,x)>0$ for all $x\in\Rm^n$ and $t>0$},
\ee
as long as the initial conditions satisfy
\be\label{feb106}
\hbox{$u_0(x):=u(0,x)\ge 0$ and $v_0(x):=v(0,x)\ge 0$ for all $x\in\Rm^n$,}
\ee
and neither $u_0(x)\equiv 0$ nor $v_0(x)\equiv 0$. 
We will always assume that (\ref{feb102}) and (\ref{feb106}) hold
 and will only consider non-negative solutions to (\ref{oct1702}).  

In addition,  the initial conditions $u_0(x)$ and $v_0(x)$ are assumed to be uniformly bounded: there exists $K>0$ so that 
at $t=0$ we have
\be\label{oct1704}
\bal
0\le u_0(x)\le K,~~0\le v_0(x)\le K,~~\hbox{ for all $x\in\Rm^n$.}
\enbal
\ee
The last assumption on the initial conditions is that there is  a uniformly positive lower bound on the sum $u_0(x)+v_0(x)$:
\be\label{oct2428}
0<K^{-1}\le u_0(x)+v_0(x)\le K,~~\hbox{ for all $x\in\Rm^n$.}
\ee
For example, assumption (\ref{oct2428}) is satisfied if the fuel is present everywhere initially: $u_0(x)\ge K$ for all $x\in\Rm^n$.
It also holds if the initial conditions
are front-like, in the sense that 
\be\label{23feb1302}
u_0(x)+v_0(x)=K,
\ee 
where $K>0$ is a non-negative constant, and there exists a direction $e\in\Sm^{n-1}$ so that $u_0(y+re)\to 0$ as $r\to+\infty$ and $v_0(y+re)\to 0$ as $r\to-\infty$,
for each $y\in\Rm^n$ fixed.  

We will make two extra assumptions on the nonlinearity~$g(v)$, in addition to (\ref{feb102}):~$g(0)=0$ and there exist $c_1>0$, $c_2>0$, $\rho\in[0,1)$ and $Z>0$ so that 
\begin{equation}\label{oct1708}
\begin{aligned}
c_1v\le g(v)\le  c_2 \exp (Zv^\rho ) \text{ for all $v\ge 0$.} 
\end{aligned}
\end{equation}
To avoid minor extra technicalities, we also suppose that the function $g(v)$ is smooth. It follows
that there exists a constant~$c_3>0$ so that
\be\label{oct2704}
c_1v\le g(v)\le c_3v,~~\hbox{ for all $0\le v\le 1$.}
\ee 

The positivity of $g(v)$ guarantees a uniform upper bound
for $u(t,x)$. More, precisely,  
the maximum principle for the standard heat equation implies that
\be\label{feb704}
0\le u(t,x)\le K, \hbox{ for all $t\ge 0$ and $x\in\Rm^n$},
\ee
as long as (\ref{oct1704}) holds at $t=0$.  

Let us comment that one can consider systems with more general nonlinearities of the form
\be \label{feb202}
\begin{split}
\partial_t u = \nu\Delta u - f_1(u,v),\\
\partial_t v = \kappa \Delta v +f_2(u,v),
\end{split}
\ee
with appropriate assumptions on the nonlinearities $f_{1,2}(u,v)$ but in this paper we prefer to keep the setting as simple as possible.

\subsubsection*{Global bounds on the solutions: bounded domains}

Let us first recall some of the existing results on global existence and bounds 
for the solutions to the thermo-diffusive systems and closely related models 
in bounded domains.
The literature is truly vast, as they arise in a wide range of applied sciences, from chemistry to atmospheric sciences. 

The question of global existence and uniform bound for such models was 
originally posed by R.~Martin in the 1970's.  
An important  step was the work of Hollis, Martin and Pierre \cite{HMP1987} 
on a system of the form (\ref{feb202}) with at most 
polynomially growing nonlinearities, typically with 
the assumptions that 
(i) the reactions are balanced, in the sense that 
\begin{equation}
\label{mar203}
-M\leq f_1(u,v)+f_2(u,v)\leq M,
\end{equation}
and that (ii) an {\it a priori} estimate, such as  (\ref{feb704}), holds for 
one of the unknowns, say, for $u(t,x)$.  
To give a glimpse of the needed ingredients, let us briefly 
review a simplified proof of the main results of  \cite{HMP1987} in the special case 
$M=0$ in (\ref{mar203}), and with the Dirichlet boundary conditions on $\partial\Omega$:
\begin{equation}
\label{mar205}
u(t,x)=v(t,x)=0\ \ \hbox{on $\partial\Omega$}.
\end{equation} 
The key observation is 
that an auxiliary unknown
\[
W(t,x)=\int_0^t\bigl(u(s,x)+\frac\nu\kappa v(s,x)\bigl)~ds,
\]
satisfies a forced heat  equation
\begin{equation}
\label{mar204}
\partial_tW-\nu\Delta W=
\frac\nu\kappa\bigl(u_0(x)+v_0(x)\bigl)+(1-\frac\nu\kappa)u(t,x),
\end{equation}
with the Dirichlet boundary condition $W(t,x)=0$ on $\partial\Omega$. 
As the right side of (\ref{mar204}) is bounded 
by assumption~(ii) above, and the domain~$\Omega$ is bounded, the function $W(t,x)$ is 
uniformly bounded, by  a comparison to a multiple of the time-independent
solution to 
\be
-\nu\Delta\bar W=1,
\ee
with $\bar W=0$ on $\partial\Omega$. 
The standard parabolic regularity results in bounded domains~\cite{Lieberman}
imply then that~$\Vert v(t,\cdot)\Vert_{L^p(\Omega)}$ is   uniformly bounded for 
all $p\in(1,+\infty)$. Because $f_2(u,v)$ has at most polynomial growth in $v$, 
while $u$ is bounded 
by assumption~(ii), we deduce that
$\Vert f_2\bigl(u(t,\cdot),v(t,\cdot))\bigl\Vert_{L^p(\Omega)}$ is also uniformly bounded for all $p\in(1,+\infty)$. Another application of parabolic regularity 
gives the uniform boundedness of $v$, finishing the proof.  

In \cite{HMP1987}, the bounds are derived from a more sophisticated duality method that 
has had applications to more general systems than  (\ref{feb202}). Surprisingly,  
the 
assumption (\ref{mar203}) cannot be much relaxed, as 
Pierre and Schmitt \cite{PS1997} constructed examples of $f_1$ and $f_2$ 
where 
the left side bound 
in~(\ref{feb202}) does not hold, and where solutions blow up in a finite time.

When $f_1$ and/or $f_2$ grow faster than any polynomial, less is known. Global existence of weak solutions for systems of the form (\ref{feb202})   is proved in Pierre \cite{P1987} under the additional conditions that one of the $f_i$'s grows linearly, while no assumption is made on the other. Whether the solutions are classical is not known. The question is also still open for an arbitrary $g(v)$ in the case of our original system (\ref{oct1702}). We refer to \cite{P2010} for an extensive literature review, especially for numerous results on bounded domains,
with either Dirichlet or Neumann boundary conditions.

Another line of investigation concern systems modeling reversible reactions. With no intention of being general, let us provide an example of such systems. For a chemical process involving four reactants $A_i, 1\leq i\leq 4$ with mass fractions $Y_i$, the progress of the reaction
$$
A_1+A_2\rightleftharpoons A_3+A_4
$$
can be described by the system
\begin{equation}
\label{mar207}
\partial_tY_i-d_i\Delta Y_i=\varepsilon_i(Y_1Y_2-Y_3Y_4),
\end{equation}
with $\varepsilon_i>0$ if $i=1,2$ and $<0$ if $i=3,4$. Interestingly,  
such models may be derived from reactive Boltzmann equations, as shown by
Bisi and Desvillettes \cite{BD2006}, which gives a natural way to 
introduce an entropy. This has opened fruitful research directions,  such as global existence questions in Desvillettes et al. \cite{DFPV2007}, trend to equilibrium 
in Devillettes, Fellner and Tang~\cite{DFT2017}, 
and systems coupling inside to surface dynamics in Fellner, Latos and Tang  \cite{FLT2018}.

\subsubsection*{The case of the whole space}

The global existence proof for solutions to (\ref{oct1702}) with polynomially growing 
nonlinearities could be
reproduced essentially 
verbatim in the whole space. However, the lack of a spectral gap estimate on the Laplacian
makes the uniform bounds on the whole space a much 
more challenging problem. In addition, there is, 
potentially, an infinite supply of fuel, unlike in a bounded domain.
This, in 
principle, could lead to an unbounded growth of the temperature. On the other hand, 
there is an elegant elementary 
observation from the paper~\cite{MP1992} by Martin and
Pierre that contains many beautiful insights into reaction-diffusion systems.
\begin{prop}\label{prop-pierre} 
Consider the system (\ref{oct1702}) posed in the whole space $\Rm^n$ under the assumptions~(\ref{feb102}) and (\ref{oct1704}). 
If~$\nu\ge\kappa$ then the solution to~(\ref{oct1702}) exists for all $t\ge 0$ and there exists a constant $C_0$ such that 
\be
\hbox{$0\le u(t,x)\le C_0$ and
$0\le v(t,x)\le C_0$ for all $t>0$ and $x\in\Rm^n$.}
\ee
\end{prop}
The proof is very short and uses two simple but  crucial ideas: first, $u(t,x)$ is
bounded from above and remains non-negative by the comparison principle. Second, if $\nu\ge\kappa$
then the heat kernels have the comparison
\be
\nu^{n/2}G_\nu(t,x)\ge \kappa^{n/2}G_\kappa(t,x),~~\hbox{for all $t>0$ and $x\in\Rm^n$.}
\ee
These two facts together with the Duhamel formula imply the result more or less
immediately. We recall the details in Section~\ref{sec:compare}. 
Here, we used the notation  
\be\label{apr202}
G_\kappa(t,x)=\farc{1}{(4\kappa\pi t)^{n/2}}e^{-|x|^2/(4\kappa t)},
\ee
for the standard heat kernel. 
We will use the convention that $G_\kappa(t,x)=0$ for $t\le 0$.

Thus, the problem of uniform bounds for the solutions to (\ref{oct1702}) is trivial in the case $\nu\ge\kappa$ which seemingly
can allow the concentration to arrive very fast into the places where the temperature is high, promoting its growth. It is, strangely, quite difficult 
in the physically supposedly simpler situation when~$\nu\le\kappa$. We will solely focus on that situation in the rest of this paper.

A localized Lyapunov functional  was adopted by Collet and Xin in~\cite{CX1996}
to obtain the global existence of solutions to (\ref{oct1702}) in $\Rm^n$
with bounded non-negative initial conditions 
for polynomial nonlinearities of the form $g(v)=v^m$. They established an  $L^\infty$-bound
of the form 
\be
\|v(t,\cdot)\|_{L^\infty}\le C\log\log t,~~\hbox{as $t\to+\infty$.} 
\ee  
 The authors of \cite{CX1996} considered nonlinearities of the 
form $uv^m$, but their argument extends to $u^m e^v$ with appropriate modifications. The Lyapunov functional method was used to prove a 
global bound in a more elaborate system with more variables with a polynomial nonlinearity as well in \cite{K2001}.  

Existence of the global in time solution for the exponential 
nonlinearities  $g(v)=\exp v$
in the same setting  was proved in~\cite{HLV} by Herrero, Lacey and Velasquez. 
They do not state explicitly any bound on the growth of the solution,
but, to the best of our understanding, the proof in~\cite{HLV} produces an estimate
of the type 
\be
\|v(t,\cdot)\|_{L^\infty}\le C\log t,~~\hbox{as $t\to+\infty$.} 
\ee  
Chen and Qi considered in~\cite{CQ2008} the Fisher-KPP thermo-diffusive system, with $g(v)=v$ and in one dimension:
\be  
\begin{split}
&\partial_t u = \nu \partial_x^2u - uv,\\
&\partial_t v = \kappa \partial_x^2v + uv,
\end{split}
\ee
with the initial condition $u_0(x)\equiv 1$. 
They were able to show not only that the solution is uniformly bounded in time but also that the front has a finite width
and is located at the position 
\be
m(t)=2c_*t-\farc{3}{2\lambda_*}\log t+x_0,~~\hbox{as $x\to+\infty$,}  
\ee
with $c_*=2\sqrt{\kappa}$, $\lambda_*=1/\sqrt{\kappa}$. 
This is exactly the Bramson asymptotics for the location of the front of the solution to the single Fisher-KPP equation
\be
v_t=\kappa v_{xx}+v-v^2,
\ee
with a compactly supported initial condition $v_0(x)$, originally proved in~\cite{Bramson1,Bramson2}. Bramson's results were recently extensively
revisited both  in the PDE and probabilistic literature -- we refer 
to~\cite{Graham} for the finest asymptotics and further references. 
While the proof in~\cite{CQ2008} relied
heavily on the probabilistic results by Bramson for the single equation, the above recent results would allow to adapt their proof to purely analytic 
arguments.  One key observation of~\cite{CQ2008} will be used extensively in the present paper. 

We should also mention that related bounds for thermo-diffusive
systems and traveling waves in cylinders in systems with advection,  with extra dissipative mechanisms, 
have been considered, for instance, in~\cite{BHKR2005,BN1992,D2021,HR2005,HR2010}. However, the extra dissipation played a key role in most arguments 
in these references.


\subsubsection*{The main result}

The main result of the paper is the following uniform bound for the solutions to the thermo-diffusive system. 
\begin{thm}\label{thm:main} 
Let $u(t,x),v(t,x)$ be a solution to \eqref{oct1702} with the initial conditions $u(0,x)=u_0(x)$ and~$v(0,x)=v_0(x)$ 
that satisfy~(\ref{oct1704}) and (\ref{oct2428}), with some $K>0$. Suppose that $g(v)$
satisfies the assumption \eqref{oct1708} and is smooth. 
Then, $u(t,x)$ and $v(t,x)$ are globally bounded in time: there exists~$C>0$ that depends
on the constants in the above assumptions, so that   
\begin{equation}\label{jul2218}
0 \le u (t,x)\le C, ~~0 \le v(t,x) \le C,
\end{equation}
for all $t>0$ and $x\in\Rm^n$. 
\end{thm}

\subsubsection*{Outline of the proof and organization of the paper}

Let us now briefly explain the main ingredients in the proof of Theorem~\ref{thm:main}. It is not difficult to see that, 
given the uniform bound (\ref{feb704}) on the fuel concentration, 
the first essential issue is the following. Let  $f(t,x)$ and $h(t,x)$ be the solutions to the standard heat equations 
\be\label{mar2206}
\partial_tf=\nu\Delta f+Z(t,x),
\ee
and
\be\label{mar2208}
\partial_th=\kappa\Delta h+Z(t,x),
\ee
with the initial conditions $g(0,x)=h(0,x)=0$. Note that (\ref{mar2206}) and (\ref{mar2208}) have two different diffusivities~$\kappa\neq\nu$ but the forcing
term $Z(t,x)$ in the two equations is the same. Suppose that we know that~$Z(t,x)\ge 0$ for all $t\ge 0$ and $x\in\Rm^n$ and that $f(t,x)$ satisfies
a uniform bound
\be\label{mar2210}
0\le f(t,x)\le K,~~\hbox{for all $t\ge 0$ and $x\in\Rm^n$.}
\ee
The question is what kind of bound on the function $h(t,x)$ one can infer from (\ref{mar2210}).  The main observation of~\cite{HLV} 
is that one can not deduce a uniform bound on $h(t,x)$ but it does obey a parabolic BMO bound: see Proposition~\ref{lem-jul2202} in Section~\ref{sec:k>1}.
The basic reason for this estimate is that $h$ is related to~$f$ by a singular integral type operator that maps $L^\infty$ to a BMO space. 
The proof of Proposition~\ref{lem-jul2202} below,
which is the first crucial step in the proof of Theorem~\ref{thm:main},  is different from that in~\cite{HLV} and seems to us to be
easier to follow.  This parabolic BMO bound does not use essentially any information on the forcing~$Z(t,x)$ except for (\ref{mar2210}).  

The next step is to use the parabolic BMO bound on the solution, together with 
the John-Nirenberg inequality, the upper bound on $g(v)$ in (\ref{oct1708})
and 
standard parabolic regularity estimates, to obtain an estimate that controls 
the pointwise values of the fuel temperature $v(t,x)$ in a parabolic cylinder
by its averages over a slightly larger parabolic cylinder. The precise statement is found in Corollary~\ref{cor-oct2506}
in Section~\ref{sec:av-point}. This step only requires an upper exponential bound 
on the reaction rate~$g(v)$. 
Thus, the proof of Theorem~\ref{thm:main} is reduced to bounding the space time averages 
\be\label{feb1402}
(v)_Q=\farc{1}{|Q|}\int_Q v(t,x)dxdt,
\ee
where $Q$ is a parabolic cylinder (see Section~\ref{sec:k>1} for the precise definitions). 

In order to bound the averages $(v)_Q$, we first
obtain some properties of the solutions
to (\ref{oct1702}) that are consequences of 
the hair trigger effect 
for the thermo-diffusive system~(\ref{oct1702}), described in
Section~\ref{sec:compare}. The term is borrowed from the seminal paper of 
Aronson and Weinberger \cite{Aronson-W}. It accounts for the fact that the rest state $v\equiv0$ in the Fisher-KPP equation is unstable with respect 
to any nontrivial small perturbation. An ingenious observation of~\cite{CQ2008}, paralleling that in Proposition~\ref{prop-pierre} above, 
is that the fuel temperature~$v(t,x)$ is a super-solution to
a differential inequality of the form 
\be\label{feb1404}
v_t\ge\kappa\Delta v+(d_0-d_1v)g(v),
\ee
with some constants $d_0>0$ and $d_1>0$ that depend on the initial conditions $u_0$ and $v_0$. It is here that one uses the lower bound
(\ref{oct2428}) on the sum $u_0+v_0$. In~\cite{CQ2008}, that considered the special case~$g(v)=v$ in one spatial dimension, this inequality, together with Bramson's 
estimates on the front location for the single Fisher-KPP equation in one dimension, essentially finished the proof because it gives tight 
bounds on the region where anything non-trivial can happen.   
In the general case $g(v)\neq v$, the situation is very different 
but (\ref{feb1404}) is still crucial. 
The reason is that (\ref{feb1404}) implies that, as soon as the local average~$(v)_Q$ is not too small, the solution spreads at a speed that may be very small but positive. This is the hair trigger effect
and the lower bound in assumption (\ref{oct1708}) for $v>0$ but small is critical here. This effect gives us three estimates: first, a forward in time lower bound
on $v(t,x)$ described in Proposition~\ref{prop-oct2402}, second, a forward in time upper bound on $u(t,x)$ in Corollary~\ref{cor-oct2402} and, third,
an absolutely crucial backward in time upper bound on $v(t,x)$ in a linearly growing region, presented in Corollary~\ref{cor-oct2404}, 

The last step is to obtain an upper bound on the spatial averages $(v)_Q$ of the fuel temperature. This is done in Section~\ref{sec:av-bound}. 
First, we  deduce in Section~\ref{sec:holder-reduce} from some algebraic manipulations that an upper bound on the spatial averages $(v)_Q$ of the fuel temperature
is equivalent to a uniform in time H\"older bound on the solution $J(t,x)$ to the forced heat equation
\be\label{23feb1410}
J_t-\kappa\Delta J=\bar u(t,x)-u(t,x),
\ee
Here, $\bar u(t,x)$ is the solution to the heat equation
\be
\bar u_t=\nu\Delta \bar u,~~\bar u(0,x)=u_0(x).
\ee
The proof of the H\"older bound on $J(t,x)$  relies on the consequences of the
hair trigger effect proved in Section~\ref{sec:compare}. They 
allow us to control the time intervals on which~$u(s,y)$,  
may behave non-trivially, for~$(s,y)$ in a parabolic cylinder
centered at $(t,x)$.  Roughly, $u(s,y)$ obeys a nice H\"older bound 
at what we call ``the early times", before reaction really starts, because $v(s,y)$ 
has to be small by the backward-in-time bound in Corollary~\ref{cor-oct2404}.
On the intermediate time scales, when reaction actually happens, we can say nothing about $u(s,y)$ apart from the fact that it is bounded. However, the
intermediate times can not be too long by the forward-in-time lower bound 
on $v(s,y)$ in Proposition~\ref{prop-oct2402}. Therefore, the contribution of this 
time interval to the H\'older constant of $J(t,x)$ is bounded.
Finally, at the late times, after the reaction has passed, the fuel
concentration $u(s,y)$ has to be small by the forward-in-time upper bound  
in 
in Corollary~\ref{cor-oct2402}. All together this leads to a H\"older bound on~$J(t,x)$. 
 This is explained in detail in Section~\ref{sec:parabola}.

{\bf Notation.} 
We use the notation $C$ for constants that may change from line to line but do not depend
on the time $t>0$.

{\bf Acknowledgement.}  We are deeply indebted to the anonymous referee of an earlier version of this manuscript for pointing out a serious 
error. JL was supported by a Royal Society University Research Fellowship 
(URF$\backslash$R1$\backslash$191492). JMR acknowledges a one 
month invitation from the Stanford Mathematics Department, at the occasion of which 
important parts of the work could be carried out. 
LR was partially supported by NSF grants DMS-1910023 and DMS-2205497, and ONR grant N00014-22-1-2174.

{\bf Data availability statement.} Data sharing is
not applicable to this article as no datasets were generated or analyzed during the current study.

\section{A parabolic BMO bound on the solution}\label{sec:k>1}

As in~\cite{HLV}, the first step in the proof of Theorem~\ref{thm:main} is a global 
in time parabolic BMO bound on the solution. Before we state the result, let us set some notation
and recall the 
definition of the parabolic~BMO semi-norm. A parabolic cylinder of size $R>0$ centered at a point~$(t,x)\in\Rm\times\Rm^n$ is a 
set of the form 
\be\label{jul2504}
Q = Q_R (t,x) = \{(s,y):~|x-y| < R, ~|t-s| < R^2 \}\subset\Rm\times\Rm^n.
\ee
Later, we will use parabolic cylinders containing past times only. We will denote such cylinders by 
\be \label{jan3101}
Q^- = Q_R ^{-} (t,x) = \{ (s,y) :~|x-y|<R, ~t-R^2 < s < t \} \subset \Rm \times \Rm^n.
\ee
For $f \in L_{loc} ^1 (\Rm\times\mathbb{R}^n)$,  we denote the average of $f$ over a measurable set $Q$ by
\begin{equation*}
(f)_Q = \frac{1}{|Q| } \int_Q f(x, t) dxdt.
\end{equation*} 
Here, $|Q|$ is the $(n+1)$-dimensional Lebesgue measure of $Q$.
The parabolic BMO semi-norm is defined~by 
\begin{equation*}
\|f \|_{\pbmo} = \sup \Big \{ \frac{1}{|Q|} \int_Q | f - (f)_Q | dxdt  \Big\},
\end{equation*}
with the supremum is taken over all parabolic 
cylinders $Q \subset \mathbb{R}^n \times \mathbb{R}_{+}.$ 
Here is the main result of this section. 
\begin{prop}\label{lem-jul2202}
Let $u(t,x)$, $v(t,x)$ be the solution to the system \eqref{oct1702}  
with the initial conditions~$u(0,x)=u_0(x)$ and~$v(0,x)=v_0(x)$ 
that satisfy~(\ref{oct1704}). Assume, in addition, that $g(v)\ge 0$ for
all $v\ge 0$. 
There exists a constant $C_1$ that depends on the constant $K$   in (\ref{oct1704}), 
and also on the dimension $n$ and the diffusivities $\kappa>0$ and $\nu>0$ but not on the function
$g(v)$ such that 
\be\label{jul2302}
\|v\|_{\pbmo(\Rm_+\times\Rm^n)}\le C_1.
\ee
\end{prop}
Let us emphasize that the proof of Proposition~\ref{lem-jul2202} uses neither the upper sub-exponential
bound nor the lower bound in (\ref{oct1708})  on the reaction rates. Its conclusion holds 
as long as the solution to~(\ref{oct1702}) exists and the production rate $g(v)$ is non-negative for $v\ge 0$.

\subsection{The proof of Proposition~\ref{lem-jul2202}}

The idea behing the proof of Proposition~\ref{lem-jul2202} is very simple.
We may use the Duhamel formula
to represent the solution to (\ref{oct1702})  as 
\be\label{jul2212}
u (t,x)=e^{t \nu \Delta}u_{0}(x)-F(t,x),~~v(t,x)=e^{t\kappa \Delta} v_{0}(x)+H(t,x).
\ee
Here, the functions $F$ and $H$ are the solutions to the heat equations
\be\label{jul2214}
\bal
&\partial_t F =\nu \Delta F +ug(v),~~t>0,\\
&F(0,x)\equiv 0,
\enbal
\ee
and 
\be\label{feb1114}
\bal
&\partial_t H =\kappa \Delta H + ug(v),~~t>0,\\
&H(0,x)\equiv 0.
\enbal
\ee
The maximum principle for the standard heat equation together with (\ref{oct1704}) implies that 
\be\label{jul2220}
0 \le \bar u(t,x):=e^{t\nu \Delta} u_0(x) \le K,~~0 \le \bar v(t,x):=e^{t\kappa \Delta} v_0(x) \le K,~~\hbox{ for all $t>0$ and $x\in\Rm^n$.} 
\ee
Hence, a parabolic BMO bound on $v(t,x)$ would follow from such bound on $H(t,x)$. 

Moreover  the comparison principle for the heat equation  implies that $0\le u(t,x)\le \bar u(t,x)$, and thus 
\be\label{jul2221}
0 \le F(t,x) \le K, ~\hbox{ for all $t\ge 0$ and $x\in\Rm^n$.}
\ee
 Note that   (\ref{jul2214}) and (\ref{feb1114}) are heat equations with identical forcing but different diffusivities. 
One may then ask: does the uniform bound on $F(t,x)$ in (\ref{jul2221})
imply some bound on the solution $H(t,x)$ to (\ref{feb1114}) even though the diffusivities in the two equations
are different? The physical issue here is the following: suppose, for instance,
that the diffusivity $\kappa$ in the equation (\ref{feb1114})  for $H$ is larger than the diffusivity $\nu$ in the equation  (\ref{jul2214})
for $F$. On one hand, the larger diffusivity should help the
solution "cool" faster, so that if $F$ is bounded, one may expect that the function $H$ should be bounded as well.
On the other, the larger diffusivity in the equation for $H$ may take the underlying Brownian particles in the Feynman-Kac representation of the solution to an inhomogeneous heat 
equation 
\be\label{jul2502}
H(t,x)=\Em_x\int_0^t  U(t-s,x+X_s)g(V(t-s,x+X_s))ds,
\ee
to the regions where the forcing~$g(v(s,x))$ is large. This would
promote the growth of $H$ that is faster than the growth of $F$. Here, $X_s$ is the Brownian motion with diffusivity $\kappa$ starting at the position~$x$ at the time~$s=0$.
This non-trivial competition between cooling by faster spreading and heating by potentially visiting the hot spots is crucial in understanding the eventual behavior of the solution, and our understanding is far from complete. In this paper, we show that 
in the setting of thermo-diffusive systems  satisfying assumption (\ref{oct1708}) 
the heating is eventually controlled by the consumption of the fuels by the reaction
but this control is far from being simple.

Thus, our goal is to analyze
\be\label{jul2224}
H (t,x) = \int_0 ^t e^{\kappa(t-s) \Delta} p(s,\cdot) (x)ds=\int_0^t \int G_{\kappa}(t-s,x-y)p(s,y)dyds
\ee
when we already know that 
\be \label{sep930_5}
0\le F(t,x) = \int_0 ^t e^{\nu(t-s) \Delta} p(s, \cdot) (x) ds=\int_0^t \int G_{\nu}(t-s,x-y)p(s,y)dyds \le K_1.
\ee
We assume here that the function $p(s,y)$ is non-negative but essentially nothing else.

The key observation is that $H$ can be written as a space-time singular integral of
the bounded function~$F$.
This is seen as as follows. The difference $\Phi = H - F$ satisfies
\be\label{jul2541}
\bal
\partial_t \Phi = \kappa \Delta H - \nu\Delta F = 
\kappa \Delta (H-F) + (\kappa-\nu) \Delta F = \kappa\Delta \Phi + (\kappa-\nu) \Delta F,  
\enbal
\ee
with the initial condition $\Phi(0,x)=0$.   Therefore, symbolically, we can write 
\be\label{jul2238}
\Phi= (\kappa-\nu) \cT_{\kappa} F,
\ee
and
\be\label{jul2327}
H= F + (\kappa-\nu) \cT_\kappa F.
\ee
Here, we have set
\be\label{jul2239}
\cT_\kappa F:=  \Delta (\partial_t - \kappa \Delta)^{-1} F.
\ee
The  main ingredient in the proof of Proposition~\ref{lem-jul2202} is the following $\pbmo$ bound.
\begin{lem} \label{jul151}
There exists a constant $C_{\cT}$, which depends on the diffiusivity $\kappa>0$ and the dimension~$n$  such that for all $F\in L^\infty(\Rm_+\times\Rm^n)$,
we have 
\begin{equation}\label{jan1102}
\|\cT_\kappa F\|_{\pbmo} \le C_\cT \|F\|_{L^\infty (\Rm_+\times\mathbb{R}^n)}.
\end{equation}
\end{lem}
This lemma corresponds to Lemma 3.3 of \cite{HLV}.
The proof presented here is different and is based on the estimates on the singular 
kernel of the operator $\cT_\kappa$. Note that Lemma~\ref{jul151}
gives immediately the parabolic~BMO bound (\ref{jul2302}) in Proposition~\ref{lem-jul2202}, as we have,
from (\ref{jul2212})  and (\ref{jul2327}):
\be\label{jul2325}
\bal
0 \le v(t,x)&= \bar v(t,x)+H(t,x) =\bar v(t,x)+F(t,x)+(\kappa - \nu) \cT_{\kappa} F(t,x) .
\enbal
\ee
The function $F(t,x)$ is uniformly bounded by (\ref{jul2221}), 
hence also bounded in the parabolic BMO space,
and the term $\cT_{\kappa}F $ is bounded in the parabolic BMO space by Lemma~\ref{jul151}.
This gives (\ref{jul2302}).  Thus, to finish the proof of
Proposition~\ref{lem-jul2202}, we only need to prove Lemma~\ref{jul151}. 

\subsection{The proof of Lemma~\ref{jul151}}

Recall that for $F \in L^\infty (\mathbb{R}^n \times \mathbb{R}_{+})$, the function $\Psi = {\cT_\kappa}F$ is the solution to
\be\label{jul2330}
\bal
&\partial_t \Psi  = \kappa\Delta \Psi + \Delta F, \\
&\Psi (x, 0) = 0.
\enbal
\end{equation} 
Extending $F$ by 0 to negative times, the action of $\cT_\kappa$  can be written as a singular integral operator 
\be 
\cT_\kappa F (t,x) = \int_{\mathbb{R}^n \times \mathbb{R}} K(t-s,x-y) F(s,y) dyds,
\ee
with the kernel
\be\label{jul2516}
K(t,x) = \Delta G_{\kappa}(t,x) = 
\frac{1}{(4\pi \kappa t )^{{n}/{2} }} e^{-{|x|^2}/{(4\kappa t)} } \Big( \frac{1}{4\kappa^{{2}}}\frac{|x|^2}{{t^2} }- \frac{n}{2\kappa t} \Big),~~ t>0 , 
\ee
and the convention that $K(t,x)=0$ for $t\le 0$. 

We need to show that there exists a constant $C_\cT$ such that for any parabolic cylinder
\be\label{jul2333}
Q = Q_R (x_0, t_0) = \{(t,x) \in \Rm\times\mathbb{R}^n \times :~ |x-x_0| < R,  |t-t_0 | < R^2 \},
\ee
and any function $F\in L^\infty(\Rm\times\Rm^n)$, we have 
\begin{equation}\label{jul2331}
\frac{1}{|Q|} \int_Q |\cT_\kappa F - (\cT_\kappa F)_Q | dxdt \le C_\cT  \| F \|_{L^\infty }.
\end{equation}
Note that for any constant $a$, we can write
\[
a-(\cT_\kappa F)_Q=\farc{1}{Q}\int_Q(a-\cT_\kappa F)dxdt,
\]
so that
\begin{equation}\label{jul2332}
\frac{1}{|Q|} \int_Q |\cT_\kappa F - (\cT_\kappa F)_Q | dxdt \le \frac{1}{|Q|} \int_Q |\cT_\kappa F - a| dxdt + |a - (\cT_\kappa F)_Q | \le \frac{2}{|Q|} \int_Q |\cT_\kappa F - a| dxdt.
\end{equation}
Therefore, it suffices to show that there exists a constant $K_\cT$,  such that for any parabolic cylinder~$Q$ and 
function~$F\in L^\infty(\Rm\times\Rm^n)$
we can find a constant $a_Q\in\Rm$ so that 
\begin{equation}\label{jul2527}
\frac{1}{|Q|} \int_Q |\cT_\kappa F - a_Q | dxdt \le K_{\cT}  \| F \|_{L^\infty }.
\end{equation}
The constant $K_\cT$ should not depend on $Q$ or $F$ but the constant $a_Q$ may. 


Given a cylinder $Q=Q_R(t_0,x_0)$, we use the notation $2Q=Q_{2R}(t_0,x_0)$ and decompose 
\be\label{jul2334}
F = (F - (F)_Q) \one_{2Q} + (F -(F)_Q) \one_{(2Q) ^c } + (F)_Q =: F^{(1)} + F^{(2)} + (F)_Q.
\ee
Note that, as $(F)_Q$ is a constant, $\cT_\kappa (F)_Q = 0$, so that
\be\label{jul2340}
\cT_\kappa F=\cT_\kappa F^{(1)}+\cT_\kappa F^{(2)}.
\ee
We will treat the two terms in (\ref{jul2340}) separately. 

\subsubsection*{A bound on $\cT_\kappa F^{(1)}$}

For the first term in (\ref{jul2340}), we write
\be\label{jul2335}
\frac{1}{|Q|} \int_Q | \cT_\kappa F^{(1)} | dxdt \le \frac{1}{|Q|}  \Big( \int_{Q} |\cT_\kappa F^{(1)}|^2 dxdt \Big)^{{1}/{2} } |Q|^{{1}/{2} } 
\le \Big(\frac{1}{|Q|} \int_{\mathbb{R}^n\times \mathbb{R}} |\cT_\kappa F^{(1)}|^2 dxdt \Big)^{{1}/{2}}.
\ee
To bound the right side of (\ref{jul2335}), let us show that the operator $\cT_\kappa$ is bounded on $L^2(\Rm\times\Rm^n)$. This is a consequence
of Mikhlin's theorem~\cite{S1993} but let us give a short direct proof. Recall that  
\be\label{jul2336}
\cT_\kappa F = \Delta \tilde{\Phi},
\ee
and the function $\tilde\Phi$ satisfies
\begin{equation}\label{may23-2}
\partial_t \tilde{\Phi} = \kappa\Delta \tilde{\Phi} + F,~~ t>0,
\ee
with
$\tilde{\Phi} (x,t) = 0$ for  $t \le 0$.
Multiplying \eqref{may23-2} by $ \Delta \tilde{\Phi}$, integrating, we obtain
\[
\frac{1}{2} \int |\nabla \tilde{\Phi} | ^2 (0,x) dx - \frac{1}{2} \int |\nabla \tilde{\Phi} | ^2 (T,x) dx = 
\kappa \| \Delta \tilde{\Phi} \|_{L^2 (\mathbb{R}^n \times [0, T])} ^2 + \int_{\mathbb{R}^n \times [0, T]} F \Delta \tilde{\Phi} dxdt, 
\]
As $\nabla \tilde{\Phi} (0,x) = 0$, this gives
\be
\bal
\kappa\| \Delta \tilde{\Phi} \|_{L^2 (\mathbb{R}^n \times [0, T])} ^2 +  
\frac{1}{2} \int |\nabla \tilde{\Phi} | ^2 (T,x) dx \le \| F\|_{L^2 (\mathbb{R}^n \times [0, T] ) } \| \Delta \tilde{\Phi} \|_{L^2 (\mathbb{R}^n \times [0, T])} , 
\enbal
\ee
whence
\be
\bal
 \| \Delta \tilde{\Phi} \|_{L^2 (\mathbb{R}^n \times [0, T])} \le \frac{1}{\kappa} \| F\|_{L^2 (\mathbb{R}^n \times [0, T] ) }.
 \enbal
 \ee
 Taking into account (\ref{jul2336}) and passing to the limit $T\to+\infty$, we get an $L^2-L^2$ bound on the operator~$\cT_\kappa$:
 \be\label{jul2337}
 \|\cT_\kappa F\|_{L^2 (\mathbb{R}^n \times \mathbb{R} ) } \le \frac{1}{\kappa} \| F \|_{L^2 (\mathbb{R}^n \times \mathbb{R} ) }.
\ee
Going back to (\ref{jul2335}), we obtain
\begin{equation}\label{jul2338}
\bal
\frac{1}{|Q|} \int_Q | \cT_\kappa F^{(1)} | dxdt &\le \Big(\frac{1}{|Q|} \int_{\mathbb{R}^n\times \mathbb{R}} |\cT_\kappa F^{(1)}|^2 dxdt \Big)^{{1}/{2}} 
\le \frac{1}{\kappa} \frac{1}{|Q|^{{1}/{2} } } \Big( \int_{\mathbb{R}^n \times \mathbb{R} } |F^{(1)}|^2 dxdt  \Big)^{{1}/{2}} \\
&= \frac{1}{\kappa} \Big( \frac{1}{|Q|} \int_{2Q} |F - (F)_Q |^2 dxdt \Big)^{{1}/{2} } \le \farc{2}{\kappa}\|F\|_{L^\infty (\mathbb{R}^n \times \mathbb{R}_+ ) }.
\enbal
\end{equation}
Thus, taking $a=0$ in (\ref{jul2332}), we see that
\be\label{feb1002}
\|\cT_\kappa F^{(1)}\|_{\pbmo}\le \farc{2}{\kappa}\|F\|_{L^\infty (\mathbb{R}^n \times \mathbb{R}_+ ) }.
\ee

\subsubsection*{A bound on $\cT_\kappa F^{(2)}$}
Next, we show that if we take 
\begin{equation}\label{jul2510}
b = \int K(x_0 - y, t_0 - s) F^{(2)}(y,s) dyds,
\end{equation}
then
\be\label{jul2506}
\farc{1}{|Q|}\int_Q |\cT_\kappa F^{(2)}-b|dxdt\le C\|F\|_{L^\infty (\mathbb{R}^n \times \mathbb{R}_+ ) },
\ee
with a constant $C$ that does not depend on $Q$ or $F$. Recalling (\ref{jul2332}) and the remark below (\ref{jul2337}), we see that
this will give the desired bound
\be\label{feb1006}
\|\cT_\kappa F^{(2)}\|_{\pbmo}\le C\|F\|_{L^\infty (\mathbb{R}^n \times \mathbb{R}_+ ) }.
\ee

We first show that 
\begin{equation}\label{jul2512}
|b| \le C_1(\kappa,n) \frac{t_0}{R^2} \|F\|_{L^\infty (\mathbb{R}^n \times \mathbb{R}_+ ) }.
\end{equation}
To this end, we use (\ref{jul2516}) to write
\begin{equation}\label{jul2514}
\bal
|b|&\le \int |K(t_0-s,x_0 - y) | |F_2(s,y) | dyds\le I\cdot \| F_2 \|_{L^\infty (\mathbb{R}^n \times \mathbb{R}_+ )} ,
\enbal
\ee
with
\be\label{jul2516bis}
\bal
I=
&\int_0 ^{t_0} \int_{|x_0 - y| \ge 2R} \frac{e^{-{|x_0 - y|^2}/{(4 \kappa(t_0 - s)) } }}{(4 \pi \kappa (t_0 - s) )^{{n}/{2} } }   
\Big(\frac{1}{4\kappa}\frac{|x_0 - y|^2}{{t_0 - s}}+ \frac{n}{2\kappa} \Big) \frac{1}{\kappa (t_0 - s)}  dyds \\
&= \int_0 ^{t_0} \int_{|y| \ge 2R} \frac{e^{-{|y|^2}/{(4 \kappa s) } }}{(4 \pi \kappa s )^{{n}/{2} } }   
\Big(\frac{1}{4\kappa}\frac{|y|^2}{s}+ \frac{n}{2\kappa} \Big) \frac{1}{\kappa s}  dyds  \le C(n,\kappa)
 \int_0 ^{t_0} \int_{|y| \ge 2R} \frac{e^{-{|y|^2}/{(8 \kappa s) } }}{s^{{n}/{2}+1 } }   
   dyds 
\\
&\le  \farc{C(n,\kappa)}{R^2}
 \int_0 ^{t_0} \int_{|y| \ge 2R} \frac{|y|^2e^{-{|y|^2}/{(8 \kappa s) } }}{s^{{n}/{2}+1 } }   
   dyds \le
   \farc{C(n,\kappa)}{R^2}
 \int_0 ^{t_0} \int_{|y| \ge 2R} \frac{e^{-{|y|^2}/{(16 \kappa s) } }}{s^{{n}/{2 } } }  
   dyds \le\farc{C(n,\kappa)t_0}{R^2}.
\enbal
\ee
Together with (\ref{jul2514}), this gives (\ref{jul2512}), since
\[
\|F^{(2)} \|_{L^\infty (\mathbb{R}^n \times \mathbb{R}_+)} \le 2 \|F \|_{L^\infty (\mathbb{R}^n \times \mathbb{R} _+) }.
\]

Next we show that (\ref{jul2506}) holds. Since $F_2(t,x)$ is supported away from $(t_0,x_0)$,  we can integrate without taking the principal value,
and use the zero integral property of $K(t,x)$ to write
\begin{equation}\label{jul2520}
\bal
&\int_Q |\cT_\kappa F^{(2)} - b| dxdt = \int_Q \Big| \int_{\mathbb{R}^n \times \mathbb{R} } (K(t-s,x-y) - K(t_0-s,x_0 -y) ) F^{(2)} (s,y) dyds \Big| dxdt \\
&= \int_Q  \Big| \int_{\mathbb{R}^n \times \mathbb{R} } (K(t-s,x-y) - K(t_0-s,x_0 -y) )( F (s,y) - (F)_Q) \chi_{(2Q)^c} (s,y) dyds \Big| dxdt  \\
&\le \int_Q  \Big( \int_{\mathbb{R}^n \times \mathbb{R} \setminus 2Q} |F(s,y) - (F)_Q| |K(t-s,x-y) - K(t_0-s,x_0 - y) |dyds\Big)dxdt \\
&= \int_{\mathbb{R}^n \times \mathbb{R} \setminus 2Q} |F(s,y) - (F)_Q| \Big( \int_Q |K(t-s,x-y) - K(t_0-s,x_0-y)| dxdt\Big)dyds \\
&\le 2 \|F \|_{L^\infty (\mathbb{R}^n \times \mathbb{R} ) } \sum_{j=1} ^\infty \int_{2^{j+1} Q \setminus 2^j Q} 
\Big( \int_Q |K(t-s,x-y) - K(t_0-s,x_0 - y) | dxdt\Big)dyds. 
\enbal
\end{equation}
The kernel $K(t,x)$ satisfies the following key bounds:  
\begin{equation}\label{may24}
\bal
&|K(t,x) | \le \frac{B_2}{(\sqrt{t}+ |x|)^{n+2} }, ~~|\nabla_x K(t,x) | \le \frac{B_2}{(\sqrt{t} + |x|)^{n+3} }, ~~
|\partial_t K(x, t) | \le \frac{B_2}{\left ( \sqrt{t} + |x| \right )^{n+4} },
\enbal
\end{equation}  
with a constant  $B_2$ that depends on $n$ and $\kappa$. 
Using \eqref{may24}, we obtain
\be\label{jul2521} 
\bal
&|K(x_1, t_1) - K(x_2, t_2) | \le |K(x_1, t_1) -K(x_2, t_1)| + |K(x_2, t_1) - K(x_2, t_2) | \\
&\le |x_1 - x_2| |(\nabla_x K) (\lambda x_1 + (1-\lambda) x_2, t_1) | \\
&+ ( |K(x_2, t_1)| + |K(x_2, t_2) |)^{{1}/{2} } |t_2-t_1|^{{1}/{2} } |\partial_t K (x_2, \lambda' t_1 + (1-\lambda' )t_2) |^{{1}/{2} } \\
&\le  \frac{B_2|x_1-x_2|}{( \sqrt{t_1} + |\lambda x_1+(1-\lambda)x_2|)^{n+3}}\\ 
&+ B_2 
\Big(\farc{|t_2-t_1|^{1/2}}{(\sqrt{t_1}+|x_2|)^{(n+2)/2}}+\farc{|t_2-t_1|^{1/2}}{\sqrt{t_2}+|x_2|)^{(n+2)/2}}\Big) 
\frac{1}{( (\lambda' t _1+ (1-\lambda') t_2 )^{{1}/{2} } + |x_2|)^{(n+4)/2} },
\enbal
\ee
with some $\lambda, \lambda' \in [0,1]$. 

We will apply (\ref{jul2521}) with 
\be\label{jul2522}
\hbox{$t_1=t-s$, $t_2=t_0-s$, $x_1=x-y$, $x_2=x_0-y$,}
\ee
as appears in the last integral in the right side of (\ref{jul2520}). 
Since the point $(t_0,x_0)$ is the center of the parabolic cylinder $Q$, and the point $(t,x)$ lies  in $Q$, while the point  $(s,y)$ lies outside the parabolic cylinder $2Q$, 
there exists a constant $B_3 = B_3 (n)>1 $ such that 
\begin{equation}\label{jul2524}
\frac{1}{B_3} \le \frac{|\lambda x + (1-\lambda) x_0 - y|}{|x_0 - y|}, \Big( \frac{|t-s|}{|t_0 - s|} \Big)^{{1}/{2} }, 
\Big( \frac{|\lambda' t + (1-\lambda') t_0 -s|}{|t_0 - s|} \Big)^{{1}/{2} } \le B_3.
\end{equation}
Inserting the points given by (\ref{jul2522}) into (\ref{jul2521}) and keeping in mind (\ref{jul2524}), we get
\begin{equation}\label{jul2525}
\bal
|K(x-y, t-s) - K(x_0 - y, t_0 - s) | \le B_4(n,\kappa) \frac{ |x-x_0| + |t- t_0|^{{1}/{2} } } { \left ( |t_0 - s|^{{1}/{2} } + |x_0 - y| \right )^{n+3} }.
\enbal
\end{equation}
Therefore, for each $(s,y) \in 2^{j+1} Q \setminus 2^j Q$, we have, since $(t,x)\in Q_R(t_0,x_0)$:
\begin{equation}\label{jul2526}
\bal
&\int_Q |K(x-y, t-s) - K(x_0 - y, t_0 - s) | dxdt \le \frac{2 B_4 R|Q|}{ \left ( |t_0 - s|^{{1}/{2} } + |x_0 - y| \right )^{n+3} } 
\le \frac{B_5}{2^{(n+3) j}},
\enbal
\end{equation}
with the constant $B_5$ that depends on $n$ and $\kappa$. 
Going back to (\ref{jul2520}) and uisng (\ref{jul2526}) gives
\begin{equation}\label{jul2527bis}
\bal
&\int_Q |\cT_\kappa F^{(2)} - b| dxdt \le 2 B_5\|F \|_{L^\infty  (\mathbb{R}^n \times \mathbb{R} ) } \sum_{j=1} ^\infty  
\frac{|2^{j+1} Q \setminus 2^j Q|}{2^{(n+3) j }}\le 2 B_5 \|F \|_{L^\infty  (\mathbb{R}^n \times \mathbb{R} ) } \sum_{j=1} ^\infty  
\frac{2^{(j+1)(n+2)}}{2^{(n+3) j }}  |Q| \\
&= 2^{n+3} B_5 \|F \|_{L^\infty  (\mathbb{R}^n \times \mathbb{R} ) } |Q| = B_6 (n,\kappa)  \|F \|_{L^\infty  (\mathbb{R}^n \times \mathbb{R} ) } |Q| .
\enbal
\end{equation}
This gives (\ref{jul2506}).

With the bounds (\ref{jul2338}) and (\ref{jul2506}) in hand, we get (\ref{jul2527}) with $a_Q=b$ simply by writing
\begin{equation}\label{jul2528}
\bal
\frac{1}{|Q|}  \int_Q |\cT_\kappa F - b | dxdt&=\frac{1}{|Q|}  \int_Q |\cT_\kappa F^{(1)}+\cT_\kappa F^{(2)} - b | dxdt\! \\
&\le 
\frac{1}{|Q|}  \int_Q |\cT_\kappa F^{(1)}| dxdt +\!\frac{1}{|Q|}  \int_Q |\cT_\kappa F^{(2)} - b | dxdt \\
&\le \farc{2}{\kappa}  \|F \|_{L^\infty  (\mathbb{R}^n \times \mathbb{R} ) } +B_6 (n,\kappa)  \|F \|_{L^\infty  (\mathbb{R}^n \times \mathbb{R} ) },
\enbal
\end{equation}
finishing the proof of Lemma~\ref{jul151}.~$\Box$

\section{From bounded averages to bounded pointwise values}\label{sec:av-point}

%

In this section, we show that the point-wise values of the function 
$v(t,x)$ are bounded by its local averages. 
The first step is to show that if the average of the function $v(t,x)$ is small
in a unit parabolic cube then the higher $L^p$-norms of $v(t,x)$ in that ball are also small.
Let us recall the representation (\ref{jul2325}) for the function $v(t,x)$:
\be
v(t,x)=\bar v(t,x)+F(t,x)+(\kappa-\nu)\Psi(t,x)
\ee
Here, $\bar v(t,x)$ is the solution to the heat equation, as defined in (\ref{jul2220}), the function 
$F(t,x)$ is the solution to (\ref{jul2214}) and is uniformly bounded by (\ref{jul2221}), 
while $\Psi(t,x)={\cT}_\kappa F(t,x)$ satisfies (\ref{jul2330}):
\be\label{oct2512bis}
\bal
&\Psi_t-\kappa\Delta \Psi=\Delta F,\\
&\Psi(0,x)=0.
\enbal
\ee
We will need the 
following exponential moment bound.
\begin{prop}\label{jul152}  
For any  $Z\ge 0$ and $\rho\in(0,1)$   
there exists a constant $C_{Z,\rho}>0$ 
that depends on all of our constants  so that 
for all $t_0>1$ and $x_0\in\Rm^n$, we have  
\begin{equation}\label{jul2531}
\farc{1}{|Q_1(t_0,x_0)|}\int_{Q_1(t_0,x_0)} e^{Zv^\rho(t,x)} dxdt 
\le  C_{Z,\rho}\exp\Big(C_{Z,\rho}\big|(\Psi)_{Q_1}\big|^\rho\Big).
\end{equation}
\end{prop}
{\bf Proof.} It follows from (\ref{jul2212})
that there exists $C>0$ such that
\be\label{jul2530}
v (t,x)\le C (K_2+  H(t,x)),
~~\hbox{for all $t>0$ and $x\in\Rm^n$.}
\ee
Here, $H(t,x)$ is the (non-negative) solution to (\ref{feb1114}). 
Thus, (\ref{jul2531}) would follow if we can show that for every $Z>0$ 
and $\rho\in(0,1)$ there exists a constant $C_Z>0$ so that
for all unit parabolic cylinders~$Q_1(t_0,x_0)$ we have  
\begin{equation}\label{jul2532}
\farc{1}{|Q_1(t_0,x_0)|}\int_{Q_1(t_0,x_0)}\exp\Big\{ZH^\rho(t,x)\Big\}dxdt 
\le  C_{Z,\rho}\exp\Big(C_{Z,\rho}\big|(\Psi)_{Q_1}\big|^\rho\Big).
\end{equation}
The proof of (\ref{jul2532}) uses  the parabolic John-Nirenberg inequality 
that bounds the deviation of a~$\pbmo$ function from its local
average. 
\begin{prop}\label{lem-john-nir} (John-Nirenberg inequality) 
There exist constants $A, B>0$, such that for all parabolic cylinders~$Q=Q_R(t_0,x_0)$ and all~$\lambda>0$, we have
\begin{equation}\label{jul2536}
\big|\{ (t,x) \in Q:~|g(t,x) - (g)_Q | \ge \lambda \} | \le B \exp \Big( -\frac{A \lambda}{\|g\|_{\pbmo} }\Big) |Q|,
\end{equation}
for all $g\in\pbmo$. Here, $|Q|$ is the $(n+1)$-dimensional Lebesgue measure of $Q$.
\end{prop}
Let us use (\ref{jul2327}) to bound $H(t,x)$ as 
\be\label{feb1302}
H(t,x)\le  C+
(\kappa-\nu) \Psi(t,x).
\ee
This leads to the following estimate 
\be\label{jul2545}
\bal 
&\int_Q \exp\Big\{ZH^\rho(t,x) \Big\}dxdt\le C 
\int_Q\exp\Big\{  CZ|\Psi(t,x)-(\Psi)_Q|^\rho+CZ|(\Psi)_Q|^\rho\Big\}dxdt\\
&\le  C\exp\Big(CZ|(\Psi)_Q|^\rho \Big)
\int_Q    
\exp\Big\{CZ|\Psi(t,x)-(\Psi)_Q|^\rho\Big\}dxdt.
\enbal
\ee 
Here, we have set $Q=Q_1(t_0,x_0)$, and denoted by $C$ constants that do not depend 
the parabolic cylinder $Q$. 

We estimate the integral in the right side above as
\be\label{feb1306}
\bal
&\int_Q \exp\Big\{CZ|\Psi (t,x)\
-(\Psi)_Q|^\rho\Big\}dxdt\\
&\le 
\sum_{q =0 } ^\infty \exp \left (  
CZ(q+1)^{\rho}  \right ) 
\Big| \{ (t,x) \in Q:~ |\Psi(t,x) - (\Psi)_Q|\in [q, q+1] \} \Big|.
\enbal
\ee
We can use the John-Nirenberg inequality to bound the right side of (\ref{feb1306}) by 
\be\label{feb1308}
\bal
&\int_Q \exp\Big\{CZ|\Psi (t,x)\
-(\Psi)_Q|^\rho\Big\}dxdt\le   C 
\sum_{q=0} ^\infty \exp\big(CZq^\rho \big)  B 
\exp \Big( - \frac{Aq}{\|\Psi\|_{\pbmo} } \Big) |Q| .
\enbal
\ee
We now use the parabolic BMO bound (\ref{jan1102}) on $\Psi(t,x)$ 
to obtain from (\ref{feb1308}) that 
\begin{equation}\label{jul2552}
\bal
&\farc{1}{|Q|}\int_Q \exp\Big\{
CZ|\Psi(t,x)-(\Psi)_Q|^\rho\Big\}dxdt\le C  
\sum_{q=0} ^\infty \exp\big(CZq^\rho \big)   
\exp \Big( - \frac{Aq}{C } \Big)   \\
 &= C\sum_{q=0} ^\infty  
 \exp \Big(  CZq^\rho- \frac{A}{C }q\Big)\le C_{Z,\rho}.
\enbal
\end{equation}
Combining this with (\ref{jul2545})  gives
%
\begin{equation*}
\farc{1}{|Q|}\int_{Q} e^{ZH^\rho(t,x)} dxdt \le C\exp\Big(CZ\big|(\Psi)_Q\big|^\rho\Big)
, 
\end{equation*}
finishing the proof of (\ref{jul2532}) and thus that of Proposition~\ref{jul152}.~$\Box$

\begin{cor}
For any  $Z\ge 0$ and $\rho\in(0,1)$   
there exists a constant $C_{Z,\rho}>0$ 
that depends on all of our constants  so that 
for all $t_0>1$ and $x_0\in\Rm^n$, we have  
\begin{equation}\label{oct2521}
\farc{1}{|Q_1(t_0,x_0)|}\int_{Q_1(t_0,x_0)} e^{Zv^\rho(t,x)} dxdt 
\le  C_{Z,\rho}\exp\Big(C_{Z,\rho} (v)_{Q_1} \Big).
\end{equation}
\end{cor}
{\bf Proof.} This follows from (\ref{jul2531}) because of (\ref{jul2212}) 
and (\ref{jul2327}) and since the function $F(t,x)$ is uniformly bounded.~$\Box$

\begin{cor}\label{cor-oct2504}
For any $\eps>0$ and $p\ge 1$ there exists a constant $C_{p,\eps}>0$ 
that depends on all of our constants  so that 
for all $t_0>1$ and $x_0\in\Rm^n$, we have  
\begin{equation}\label{oct2531}
\farc{1}{|Q_1(t_0,x_0)|}\int_{Q_1(t_0,x_0)} [g(v(t,x))]^p dxdt 
\le  C_{p,\eps}[(v)_Q]^{1-\eps}\exp\big(C_{p,\eps} (v)_Q \big).
\end{equation}
\end{cor}
{\bf Proof.} Let us write, using H\"older's inequality, with some $p'>1$ and 
$q'> 1$ such 
that $1/p'+1/q'=1$:
\begin{equation}\label{oct2523}
\bal
&\farc{1}{|Q_1(t_0,x_0)|}\int_{Q_1(t_0,x_0)} [g(v(t,x))]^p dxdt\\
&\le  \farc{1}{|Q_1(t_0,x_0)|}\Big(\int_{Q_1(t_0,x_0)} v(t,x)dx\Big)^{1/p'}\Big(\int_{Q_1(t_0,x_0)} 
\Big(\frac{g^p(v(t,x))}{v^{1-1/q'}(t,x)}\Big)^{q'}dx\Big)^{1/q'}
\\
&\le  \farc{1}{|Q_1(t_0,x_0)|^{1/q}}[(v)_Q]^{1/p'}\Big(\int_{Q_1(t_0,x_0)} e^{Mv^\rho(t,x)}dx\Big)^{1/q'}.
\enbal
\ee
We used the upper bound for $g(v)$ in (\ref{oct2704}) in the last step above. 
Now, (\ref{oct2531}) follows after taking~$p'$ large enough and taking (\ref{oct2521}) into account.~$\Box$ 

Next, we recall the following version of the parabolic ABP maximum principle from~\cite{Lieberman,Picard}. 
\begin{prop}\label{prop-oct2504}
Let $h(t,x)\ge 0$ be the solution to 
\be
h_t=\Delta h+f(t,x).
\ee
Then, we have, for any $p\ge 1$, $t_0\ge 3$ and $x\in\Rm^n$: 
\be
\sup_{Q^{{ {-}}} _1(t_0,x_0)}h\le C\big(\|h\|_{L^p(Q^{ {-}}_2(t_0,x_0))}+\|f\|_{L^{n+1}(Q^{{{-}}}_2(t_0,x_0))}\big).
\ee
\end{prop}

Let us go back to the equation for $v(t,x)$:
\be
v_t=\kappa\Delta v+ug(v).
\ee
Proposition~\ref{prop-oct2504} implies that $v(t,x)$ obeys a bound
\be\label{oct2525}
\sup_{Q^{{ {-}}}_1(t_0,x_0)}v(t,x)\le C\big(\|v\|_{L^p(Q_2(t_0,x_0))}+\|ug(v)\|_{L^{n+1}(Q_2(t_0,x_0))}\big).
\ee
We deduce the following. 
\begin{cor}\label{cor-oct2506}
For any $\eps>0$ there is $C_\eps$ such that for any $t_0\ge 1$ and $x_0\in\Rm^n$ we have 
\be\label{oct2526}
\sup_{Q^{{ {-}}}_1(t_0,x_0)}v(t,x)\le C_\eps(v)_{Q_2(t_0,x_0)}^{1-\eps}\exp\big(C_{p,\eps} (v)_{Q_2(t_0,x_0)} \big). 
\ee
\end{cor}
{\bf Proof.} This is a consequence of Corollary~\ref{cor-oct2504} and Proposition~\ref{prop-oct2504}.~$\Box$ 

Thus, the problem of a global uniform bound on $v(t,x)$
is reduced to bounding the local averages of~$v$.

\section{Consequences of the hair trigger effect}\label{sec:compare}

In this section, we quantify the hair trigger effect for the thermo-diffusive system (\ref{oct1702}).
This is the claim of Proposition~\ref{prop-oct2402} that follows from the simple observation in Proposition~\ref{prop-oct2802}. 
It has two important consequences: the forward in time bound on the fuel concentration $u(t,x)$ in Corollary~\ref{cor-oct2402}
and a backward in time  upper bound on the fuel temperature in Corollary~\ref{cor-oct2404}. Both of them will play a role in Section~\ref{sec:av-bound} 
in the proof of 
the bound on the space-time averages of $v(t,x)$ in Proposition~\ref{prop-oct802}.  
 
Let us first recall that a uniform in time bound for the solutions to~\eqref{oct1702} 
with $\nu \ge \kappa$ on $\Rm^n$ can be proved by comparing the heat kernels 
with the diffusivities $\nu$ and $\kappa$, with no restrictions on the function $g(v)$ whatsoever using a beautiful
idea by Martin and Pierre~\cite{MP1992}. Namely,  we have an elementary inequality 
\be\label{oct2812}
\nu^{n/2}G_{\nu}(t,x)\ge\kappa^{n/2}G_{\kappa}(t,x),~~\hbox{ for all $t>0$ and $x\in\Rm^n$ if $\nu\ge \kappa>0$.}
\ee
Let $\bar u(t,x)$ and $\bar v(t,x)$ be as defined in (\ref{jul2220}). 
It follows from (\ref{oct2812}) that if $\nu\ge\kappa$ we can write
\be
u(t,x)=\bar u(t,x)-\int_0^t\int G_{\nu}(t-s,x-y)u(s,y)g(v(s,y))dyds,
\ee
and
\be
\bal
v(t,x)&=\bar v(t,x)+\int_0^t\int G_{\kappa}(t-s,x-y)u(s,y)g(v(s,y))dyds\\
&\le
\bar v(t,x)+\farc{\nu^{n/2}}{\kappa^{n/2}}\int_0^t\int G_{\nu}(t-s,x-y)u(s,y)g(v(s,y))dyds
\\
&\le\bar v(t,x)+\farc{\nu^{n/2}}{\kappa^{n/2}}(\bar u(t,x)-u(t,x))\le CK.
\enbal
\ee
Thus, in the case $\nu\ge\kappa$ the upper bound on $v(t,x)$ is very simple.

If $\nu\le\kappa$, the above argument gives a lower bound that played a crucial role in the analysis
in~\cite{CQ2008} in the special case $g(v)=v$ and dimension $n=1$
\be
v(t,x)\ge \bar v(t,x)+\farc{\nu^{n/2}}{\kappa^{n/2}}(\bar u(t,x)-u(t,x)),
\ee
or, equivalently,
\be
u(t,x)\ge \farc{\kappa^{n/2}}{\nu^{n/2}}\bar v(t,x)-\farc{\kappa^{n/2}}{\nu^{n/2}}v(t,x)+\bar u(t,x).
\ee
Note that if the lower bound in (\ref{oct2428}) holds and $\nu\le\kappa$ then we have
\be\label{oct2818}
\bal
\bar u(t,x)+{\farc{\kappa^{n/2}}{\nu^{n/2}}}\bar v(t,x)&={\farc{\kappa^{n/2}}{\nu^{n/2}}} G_{\kappa}(t,\cdot)\star v_0+G_{\nu}(t,\cdot)\star u_0 
\\
&\ge 
G_{\nu}(t,\cdot)\star v_0+G_{\nu}(t,\cdot)\star u_0=
G_{\nu}(t,\cdot)\star(u_0+v_0)\ge 
K^{-1}:=d_0.
\enbal
\ee
We deduce that 
\be
u(t,x)\ge d_0-d_1v,
\ee
with
\be\label{oct2816}
d_1=\farc{\kappa^{n/2}}{\nu^{n/2}}.
\ee
Let us summarize this as follows.
\begin{prop}\label{prop-oct2802}
Let $v(t,x)$ and $u(t,x)$ be the solutions to (\ref{oct1702}) with the initial conditions $u_0(x)$ and $v_0(x)$ that satisfy
(\ref{oct1704}) and  (\ref{oct2428}). 
If $\nu\le\kappa$ then  $v(t,x)$ satisfies a differential inequality
\be\label{oct1714}
v_t\ge \kappa\Delta v+(d_0-d_1v)g(v),
\ee
with $d_0$ and $d_1$ given, respectively, by (\ref{oct2818}) and (\ref{oct2816}). 
\end{prop}
We should stress that the conclusion of this proposition relies crucially on the lower bound on the sum~$u_0(x)+v_0(x)$ in (\ref{oct2428})
that makes the constant $d_0$ in (\ref{oct2818}) be positive.  

The main consequence of Proposition~\ref{prop-oct2802} is that it controls the function $v(t,x)$ 
from below on a set that grows linearly in time. Here, we rely on the lower bound in (\ref{oct1708}) for the function $g(v)$ for~$v>0$ small.  
We should mention that the upper bound in  (\ref{oct1708}) is not required here. 
%
\begin{prop}\label{prop-oct2402}
Let $v(t,x)$ and $u(t,x)$ be the solutions to (\ref{oct1702}) with the initial conditions $u_0(x)$ and $v_0(x)$ that satisfy
(\ref{oct1704}) and  (\ref{oct2428}).  Assume
also that $g(v)$ satisfies the lower bound 
\be\label{jan2312}
\hbox{$g(v)\ge c_1v$ for all $v\ge 0$,with some $c_1>0$.}
\ee 
There exist  $\gamma_0>0$, $\lambda_0>0$, $c_0>0$ and $C>0$
that depend on the constant~$K$ in~(\ref{oct1704}) and~(\ref{oct2428}), the diffusivities
$\kappa$ and $\nu$, and the constant $c_1$   so that the following holds. 
Assume that at a time~$t_0>0$ we have, for some $x_0\in\Rm^n$,  
\be\label{oct2502}
\int_{|x-x_0|\le 1}v(t_0,x)dx\ge\delta_0,
\ee
with some $\delta_0\in(0,1)$. 
Then, we have 
\be\label{oct2429}
v(t,x)\ge \min\Big(\gamma_0,C^{-1}\delta_0 e^{\lambda_0(t-t_0)}\Big),~~\hbox{ for all $t\ge t_0+1$ and $|x-x_0|\le 1$,}
\ee
and 
\be\label{nov412}
v(t,x)\ge\gamma_0,~~\hbox{for all $\disp t\ge t_1:= t_0+\farc{1}{\lambda_0}|\log\delta_0|+C$ and $|x-x_0|\le c_0(t-t_1)$.}
\ee
\end{prop} 
{\bf Proof.} {\bf Step 1. A sub-solution for $v(t,x)$.} 
First, we construct a sub-solution for $v(t,x)$. 
Observe that, due to the assumption (\ref{jan2312}), if we choose $\alpha_0=d_0/(2d_1)$, then there exists~$\beta_0>0$ so that 
\be\label{jan2304}
(d_0-d_1s)g(s)\ge \beta_0 s,~~\hbox{ for all $0\le s\le\alpha_0$.}
\ee
Let us now choose $L$ sufficiently large, so that the principal eigenvalue of the Dirichlet problem
\be\label{jan2314}
\bal
&-\kappa\Delta\varphi_L=\lambda_L\varphi_L,~~\hbox{ for $x\in B(0,L)$},\\
&\varphi_L=0,~~\hbox{ on $\partial B(0,L)$},
\enbal
\ee
satisfies 
\be\label{jan2318}
0<\lambda_L<\beta_0/2.
\ee 
We normalize the non-negative
principal eigenfunction $\varphi_L(x)$ so that 
\[
\varphi_L(0)=\sup_{x\in B(0,L)}\varphi_L(x)=1.
\]
The point is that the function $\varphi_L$ satisfies the differential inequality 
\be\label{jan2302}
\kappa\Delta\varphi_L+\beta_0\varphi_L\ge 0,~~\hbox{for $x\in B(0,L)$.}
\ee
We see from (\ref{jan2304}) and (\ref{jan2302}) that for $0<\delta\le\alpha_0$ the function
$\varphi_{\delta,L}(x)=\delta\varphi_L(x-x_0)$ obeys the differential inequality 
\be\label{jan2308}
\kappa\Delta\varphi_{\delta,L}+(d_0-d_1\varphi_{\delta,L})g(\varphi_{\delta,L})\ge 0,~~\hbox{in $B(x_0,L)$.}
\ee
Thus, by a standard argument, the function 
\[
\tilde\omega_{\delta}(x)=\delta\varphi_L(x-x_0)\one(|x-x_0|\le L)
\]
is a time-independent sub-solution to (\ref{jan1608}) on the whole space $\Rm^n$. 

Consider the function $\omega_\delta(t,x)$, the solution to 
\be\label{jan1608}
\partial_t\omega_\delta=\kappa\Delta\omega_\delta+(d_0-d_1\omega_\delta)g(\omega_\delta),~~t\ge t_0+1/2,
\ee
with the initial condition 
\be
\omega_\delta(t_0+1/2,x)=\tilde\omega_\delta(x). 
\ee
As the initial condition $\tilde\omega_\delta(x)$ is a time-independent sub-solution to (\ref{jan1608}),  
the comparison principle  shows that $\omega_\delta(t_0+1/2+h,x)- \tilde\omega_\delta(x)>0$ for any $h>0$.
In addition, the function 
\[
z(t,x)=\omega_\delta(t+h,x)-\omega_\delta(t,x)
\]
satisfies a parabolic equation of the form
\be
z_t=\kappa \Delta z(t,x)+c(t,x)z,
\ee
with a bounded function  $c(t,x)$ given by
\be
c(t,x)=\frac{(d_0-d_1\omega_\delta(t+h,x))g(\omega_\delta(t+h,x))-(d_0-d_1\omega_\delta(t,x))g(\omega_\delta(t,x))}
{\omega_\delta(t+h,x){ -}\omega_\delta(t,x)}{  \one ( \omega_\delta(t+h,x) \ne  \omega_\delta(t,x)) },
\ee
and a positive initial condition
\be
z(t_0+1/2,x)=\omega_\delta(t_0+1/2+h,x)-\tilde\omega_\delta(x)\ge 0.
\ee
The comparison principle implies that $z(t,x)>0$ for all $x\in\Rm^n$ and $t> t_0+{1}/{2}$, for any $h>0$. 
Thus, the function $\omega_\delta(t,x)$ is increasing in time. 

In addition, the function 
\be
\mu(t,x)=\delta e^{(\beta_0-\lambda_L)(t-t_0-1/2)}\varphi_L(x-x_0),
\ee
satisfies
\be
\mu_t-\kappa\Delta\mu=\beta_0\mu\le (d_0-d_1\mu)g(\mu),
~~\hbox{ in $B(x_0,L)$, as long as $\|\mu(t,\cdot)\|_{L^\infty}\le\alpha_0$,}
\ee
with the boundary condition $\mu(t,x)=0$ for $x\in\partial B(x_0,L)$.  We deduce from the comparison principle that
\be\label{jan1610}
\omega_\delta(t,x)\ge \delta e^{(\beta_0-\lambda_L)(t-t_0-1/2)}\varphi_L(x-x_0),~~\hbox{ as long as $\delta e^{(\beta_0-\lambda_L)(t-t_0)}\le C$,}
\ee
with a constant $C>0$ that depends on $\alpha_0$, $\beta_0$ and $\lambda_L$. Since the function $\omega_\delta(t,x)$ is increasing in~$t$ and the eigenfunction
$\varphi_L(x)$ is uniformly bounded from below on $B(x_0,1)$, as long as $L>2$, we deduce that there exist~$\gamma>0$ and~$C>0$ so that 
\be\label{jan2323}
\omega_\delta(t,x)\ge \min\Big(\gamma,C\delta e^{(\beta_0-\lambda_L)(t-t_0-1/2)}\Big),~~\hbox{for $t>t_0+1/2$ and $|x-x_0|\le 1$.}
\ee

{\bf Step 2. Using the sub-solution.} We now use the sub-solution $\omega_\delta(t,x)$ constructed above to bound the function $v(t,x)$ from
below. Suppose that $v(t_0,x)$ satisfies assumption (\ref{oct2502}) with some~$\delta_0\in (0,1)$.
Let us fix $L>0$ so that (\ref{jan2318}) holds. 
As $v(t,x)$ is a super-solution of the heat equation:
\be\label{jan1604}
v_t\ge\kappa\Delta v,
\ee
it follows that  there exists a constant $k_L>0$ so that the assumption (\ref{oct2502}) implies that 
\be\label{jan1606}
v(t_0+1/2,x)\ge k_L\delta_0,~~\hbox{ for all $|x-x_0|\le L$.}
\ee
We choose $\delta=\min(k_L\delta_0,\alpha_0)$. Then, at the time $t=t_0+1/2$ we have 
\be
v(t_0+1/2,x)\ge \omega_\delta(t_0+1/2,x),~~\hbox{ for all $x\in \Rm^n$.}
\ee
It follows from (\ref{oct1714}) and (\ref{jan1608}) that 
\be\label{jan2320}
v(t,x)\ge \omega_\delta(t,x),~~\hbox{ for all $t>t_0+1/2$ and $x\in \Rm^n$.}
\ee
Going back to (\ref{jan2323}), we deduce that 
\be\label{jan2321}
v(t,x)\ge \omega_\delta(t,x)\ge \min\Big(\gamma,C\delta e^{(\beta_0-\lambda_L)(t-t_0-1/2)}\Big),~~\hbox{for $t>t_0+1/2$ and $|x-x_0|\le 1$.}
\ee
As $0<\delta_0<1$, we have 
\be
\delta= \min(k_L\delta_0,\alpha_0)\ge \delta_0\min(k_L,\alpha_0)=C\delta_0.
\ee
Using this in (\ref{jan2321}) gives 
\be\label{jan2324}
v(t,x)\ge   \min\Big(\gamma,C\delta_0e^{(\beta_0-\lambda_L)(t-t_0-1/2)}\Big),~~\hbox{for $t>t_0+1/2$ and $|x-x_0|\le 1$,}
\ee
which is (\ref{oct2429}). 

A consequence of (\ref{jan2324}) is that there is $C>0$ so that at the time
\be\label{jan2326}
t_1= t_0+\farc{1}{\beta_0-\lambda_L}|\log\delta_0|+C
\ee
we have 
\be\label{jan2325}
v(t_1,x)\ge\gamma,~~\hbox{for all $|x-x_0|\le 1$.}
\ee
Now, the bound (\ref{nov412}) follows from (\ref{jan2325}) and the classical
results on spreading for solutions to the Fisher-KPP type equations for solutions that are initially strictly positive
on a unit ball~\cite{Aronson-W}.  More precisely, let $\tilde v(t,x)$ be the solution to the initial value problem
\be\label{jan2402}
\bal
&\tilde v_t=\kappa\Delta\tilde v+(d_0-d_1\tilde v)g(\tilde v),~~t\ge t_1,
\\
&\tilde v(t_1,x)=\gamma\one(|x-x_0|\le 1),
\enbal
\ee
and let $c_*$ be the minimal speed for a traveling wave solution for (\ref{jan1608}).
Then, Theorems 3.1 and~5.3 of~\cite{Aronson-W} imply that for any $0<c<c_*$ we have 
\be
\liminf_{t\to+\infty}\inf_{|y-x_0|<c(t-t_1)}\tilde v(t,y)\ge \farc{d_0}{2d_1}.
\ee
Thus, there exists a time $T_0$ that depends only on $\gamma$   so that 
\be\label{jan2404}
\tilde v(t,y)\ge \farc{d_0}{4d_1},~~\hbox{ for all $t>t_1+T_0$ and $|y-x_0|\le c(t-t_1)$.}
\ee
It follows from the comparison principle, (\ref{oct1714}), (\ref{jan2325}), (\ref{jan2402}) and (\ref{jan2404}) that
\be
v(t,y)\ge\tilde v(t,y)\ge \farc{d_0}{4d_1},~~\hbox{ for all $t>t_1+T_0$ and $|y-x_0|\le c(t-t_1)$,}
\ee
which is (\ref{nov412}). 
In particular, we can  
let the speed~$c_0$ in (\ref{nov412}) be smaller than the minimal speed for a traveling wave solution for~(\ref{jan1608}).~$\Box$
%
%
%
%

Proposition~\ref{prop-oct2402} has several important corollaries that will be crucial for us. First, we have a forward-in-time upper bound on the fuel. 
\begin{cor}\label{cor-oct2402}
Assume that the assumptions of Proposition~\ref{prop-oct2402} hold.
There exist $c_2>0$ and $\gamma_2>0$ so that   the following holds. 
Suppose that at some time~$t_0>0$ the inequality (\ref{oct2502}) holds with some~$\delta_0\in(0,1)$. 
Then, we have 
\be\label{nov402}
u(t,x)\le Ke^{-\gamma_2(t-t_1)},~~\hbox{for all $t\ge t_1=t_0+|\log\delta_0|+C$ and $|x-x_0|\le c_2(t-t_1)$.}
\ee
\end{cor}
{\bf Proof.} Informally, the reason for this estimate is that 
Proposition~\ref{prop-oct2402} implies that 
$v(t,x)\ge \gamma_0$ in the region $|x-x_0|\le c_0(t-t_1)$ that is
much bigger than the domain $|x-x_0|\le c_2(t-t_1)$ which appears in (\ref{nov402}), 
by a factor growing linearly in $t-t_1$, 
if we take, say, $c_2<c_0/2$. This gives the exponential decay on the fuel $u(t,x)$.  

To make this precise, let $\gamma_0$ be a lower bound for $v(t,x)$ in the cone 
\[
\mathcal{C}_{t_1,x_0}:=\{\vert x-x_0\vert\leq c_0\vert t-t_1\vert\},
\]
as in (\ref{nov412}), and let
\be
q_0=\inf_{v\geq\gamma_0}g(v).
\ee
Then, in the cone $\mathcal{C}_{t_1,x_0}$ the function $u(t,x)$ satisfies 
\be\label{nov404}
u_t-\nu\Delta u+q_0u\leq0.
\ee
Let us look for a super-solution to (\ref{nov404}) in the form 
\be\label{nov406}
\bar u(t,x)=Me^{-\alpha (t-t_1)}\cosh~(\beta\vert x-x_0\vert),
\ee
with the positive constants $M$, $\alpha$ and $\beta$ to be chosen. 
Note that we have, in the polar coordinates
\be\label{nov408}
\bal
\bar u_t-\nu\Delta\bar u+q_0 \bar u
=\bar u_t-\nu\bar u_{rr}-\nu\di\frac{n-1}r\bar u_r+q_0 \bar u 
=\Big(-\alpha-\nu\beta^2-(n-1)\nu\beta\di\frac{\tanh ~r}r+q_0\Big)\bar u.
\enbal
\ee
As $\tanh r$ is concave for $r>0$, we have
\[
\frac{\tanh~r}r\leq1,~~~\hbox{for all $r>0$}.
\]
Using this in (\ref{nov408}) gives 
\be\label{nov410}
\bar u_t-\nu\Delta\bar u+q_0 \bar u\geq(q_0-\alpha-\nu\beta^2-(n-1)\nu\beta)\bar u.
\ee
Therefore, we may choose $\alpha>0$ and $\beta>0$  suitably small so that  the 
right side in (\ref{nov410}) is positive.  In addition, we may choose 
$\alpha$ and $\beta$ so that we have the inclusion
\[
\tilde{\mathcal{C}}_{t_1,x_0}:=
\{\vert x-x_0\vert\leq\frac\beta\alpha\vert t-t_1\vert\}\subset \mathcal{C}_{t_1,x_0}
=\{\vert x-x_0\vert\leq c_0\vert t-t_1\vert\}.
\]
With this choice of $\alpha>0$ and $\beta>0$, we know that $u(t,x)$ is a sub-solution
and $\bar u(t,x)$ is a super-solution to (\ref{nov404}) 
in the cone $\tilde{\mathcal{C}}_{t_1,x_0}$. Moreover, as $u(t,x)\le K$ for all $t\ge 0$
and $x\in\Rm^n$, 
we can choose $M>0$ sufficiently large, so that $u(t,x)\le \bar u(t,x)$ for all 
$(t,x)$ on the boundary of $\tilde{\mathcal{C}}_{t_1,x_0}$.
It follows that then $u(t,x)\le \bar u(t,x)$ 
in the whole cone $\tilde{\mathcal{C}}_{t_1,x_0}$.
This gives the exponential decay of $u$ in any cone of the form
$\{\vert x-x_0\vert\leq c\vert t-t_1\vert\}$ with $c<\alpha/\beta$.~$\Box$

Next, we state a backward-in-time bound that shows that if the average of the 
temperature is small in a ball then the averages over unit 
balls have to be small in a growing region at earlier times. We denote by $|B_1|$ the volume of the unit ball
in $\Rm^n$. 
\begin{cor}\label{cor-oct2404}
Assume that the assumptions of Proposition~\ref{prop-oct2402} hold, 
and let
$c_0$and $\gamma_0$  be as in the conclusion of that proposition, and 
let $\eps_0 \in \left (0 , {\gamma_0}/{|B_1|} \right )$. Then there exist constants~$T_0>0$, $C_1>0$, depending on $\eps_0$ and $\gamma_3>0$
that is independent of $\eps_0$ so that the following holds. 
Assume that at some time~$t_0>T_0$ we have, for some $x_0\in\Rm^n$: 
\be\label{oct2504}
\int_{|x-x_0|\le 1}v(t_0,x)dx\le\eps_0. 
\ee
Then, we have 
\be\label{oct2432}
\int_{|y-x|\le 1}v(t,y)dy\le C\eps_0 e^{-\gamma_3(t_0-t)},
~~\hbox{for all $t\le t_0- T_0 $ and } |x-x_0|\le \farc{c_0}{2}(t_0-t)-C_1.
\ee
Moreover, we can choose 
\be
T_0 = C|\log \eps_0 | + C, ~ C_1 = C + C|\log \eps_0|
\ee
where $C$ is a constant independent of $\eps_0$.
\end{cor}
{\bf Proof.}  Let $\gamma_0$ be as in Proposition~\ref{prop-oct2402} and 
fix some $\eps_0\in(0,\gamma_0/10)$.  
Assume that~(\ref{oct2432}) is violated at some time $t'<t_0 -T_0$ and $x'\in\Rm^n$: 
\be\label{nov414}
\int_{|y-x'|\le 1}v(t',y)dy\ge C\eps_0 e^{-\gamma_3(t_0-t)}.
\ee
It follows from Proposition~\ref{prop-oct2402} that then
\be\label{nov416}
v(t,x)\ge\gamma_0,~\hbox{for all $\disp t\ge t_1':= t'+
\farc{1}{\lambda_0}|\log\eps_0| +\farc{\gamma_3}{\lambda_0}(t_0-t')+C$ 
and $|x-x'|\le c_0(t-t_1')$.}
\ee
If we choose $\gamma_3<\lambda_0$ then we have
\be
\bal
t_1':= t'+
\farc{1}{\lambda_0}|\log\eps_0| +\farc{\gamma_3}{\lambda_0}(t_0-t')+C<t_0,
\enbal
\ee
provided that $T_0>0$ is large enough: 
that is, if
\be
T_0 \ge \left ( 1 - \frac{\gamma_3}{\lambda_0 } \right )^{-1} \left ( \frac{1}{\lambda_0} |\log \eps_0 | + \frac{1}{\lambda_0} |\log C| + C \right ).
\ee
We deduce from  (\ref{oct2504}) and since $\eps_0 < {\gamma_0}/{|B_1|}$
that we
must have
\be
\bal
|x_0-x'|>c_0(t_0-t_1'){-1}=
c_0\Big(t_0-t'-\farc{1}{\lambda_0}|\log\eps_0| -\farc{\gamma_3}{\lambda_0}(t_0-t')-C\Big)
\ge \farc{c_0}{2}(t_0-t')-{C_1}.
\enbal
\ee
Now, (\ref{oct2432}) follows.~$\Box$

\section{A bound on the averages}\label{sec:av-bound} 
 
The last step in the proof of Theorem~\ref{thm:main}  is to prove a bound on the space-time averages of the function $v(t,x)$.
As we have seen in Corollary~\ref{cor-oct2506}, such bound would provide 
a global in time point-wise bound on $v(t,x)$ as well, finishing the proof of Theorem~\ref{thm:main}.   
\begin{prop}\label{prop-oct802}
Let $u(t,x)$, $v(t,x)$ satisfy the system \eqref{oct1702}  
under the assumptions~(\ref{oct1704}),~(\ref{oct2428}) and~(\ref{oct1708}).  
There exists a constant $C>0$ that depends on the constants in the above 
assumptions, and also on the dimension $n$, and diffusivities $\kappa$ and $\nu$,  
such that for any parabolic cylinder~$Q_1(t_0,x_0)$, with $x_0\in\Rm^n$
and $t_0>1$, we have
\be\label{oct810}
|(v)_{Q_1}|\le C.
\ee
\end{prop}
The rest of this section contains the proof of Proposition~\ref{prop-oct802}.

\subsection{Reduction to a H\"older in time bound}\label{sec:holder-reduce}  

Once again, let us recall the representation (\ref{jul2325}) for the function $v(t,x)$:
\be
v(t,x)=\bar v(t,x)+F(t,x)+(\kappa-\nu)\Psi(t,x)
\ee
Here, $\bar v(t,x)$ is the solution to the heat equation, as defined in (\ref{jul2220}), the function 
$F(t,x)$ is the solution to (\ref{jul2214}) and is uniformly bounded by (\ref{jul2221}), while $\Psi(t,x)$ satisfies 
\be\label{oct2512}
\bal
&\Psi_t-\kappa\Delta \Psi=\Delta F,\\
&\Psi(0,x)=0.
\enbal
\ee
Next, let $J(t,x)$ be the solution to
\be\label{feb1410}
J_t-\kappa\Delta J=F,
\ee
with the initial condition
$J(0,x)=0$. 
Then, we claim that
\be\label{oct1402}
\Psi(t,x)=\farc{1}{\kappa}(\partial_tJ-F).
\ee
 Observe that
the derivative $\partial_t J$ satisfies
\be
\partial_t (\partial_t J) - \kappa \Delta (\partial_t J) = \partial_t F, 
\ee
with the initial condition 
$(\partial_t J) (0,x) = F(0,x)$.
Hence, $(\partial_t J - F)$ satisfies 
\be\label{oct2402}
\partial_t (\partial_t J - F) - \kappa\Delta (\partial_t J - F) = \partial_t F - \partial_t F + \kappa \Delta F = \kappa\Delta F, 
\ee
with $(\partial_t J - F ) (0,x) = 0$.  This gives (\ref{oct1402}). 

Altogether we now have
\be\label{oct2510}
v(t,x)=\bar v(t,x)+F(t,x)+(\kappa-\nu)\farc{1}{\kappa}(\partial_tJ(t,x)-F(t,x))=\bar v(t,x)
+\farc{\nu}{\kappa}F(t,x)+\farc{\kappa-\nu}{\kappa}\partial_tJ(t,x).
\ee
As the functions $\bar v(t,x)$ and $F(t,x)$ are bounded, it follows that for any parabolic cube $Q$ we have
\be\label{oct2726}
(v)_{Q}\le C(1+(\partial_tJ)_Q).
\ee
Let us write
\be\label{oct2404}
\bal
(\partial_t J)_{Q_R} &=   \frac{1}{|Q_R|} 
\int_{|x-x_0| < R} \int_{|t-t_0| < R^2} \partial_t J (t,x) dt dx   \\
&=\frac{1}{R^{n+2} |Q_1|} \int_{|x-x_0|<R}\left (J(t_0 +R^2,x) - J(t_0-R^2,x)\right) dx.
\enbal
\ee
Thus the boundedness of the averages of $v(t,x)$ and the conclusion of Proposition~\ref{prop-oct802}
would follow from the following H\"older in time bound on $J(t,x)$. 
\begin{lem}\label{lem-oct2402}
For any $\alpha\in(0,1)$, there exists a constant $C_\alpha>0$ 
so that
\be\label{oct2410}
|J(t_2,x)-J(t_1,x)|\le C_\alpha|t_2-t_1|^{\alpha},
~~\hbox{ for all $t_1+1\ge t_2\ge t_1\ge 2$ and $x\in\Rm^n$.}
\ee
\end{lem}
The rest of this section contains the proof of Lemma~\ref{lem-oct2402}.

\subsection{Reduction to a bound on the spatial variations of the fuel in a "parabola"}

We now begin the proof of Lemma~\ref{lem-oct2402}. 
Assume that $t_2\geq t_1\ge 1$, and write
\begin{equation}\label{jul2616}
\bal
J(t_1,x)-J(t_2,x)=&\int_0^{t_1}\int_{\Rm^n}\Big(
\frac{e^{-\vert x-y\vert^2/4\kappa (t_1-s)}}{(4\pi\kappa (t_1-s))^{n/2}}
-\frac{e^{-\vert x-y\vert^2/4\kappa(t_2-s)}}{(4\pi\kappa(t_2-s))^{n/2}}\Big)F(s,y)dyds\\
&-\int_{t_1}^{t_2}\int_{\Rm^n}\frac{e^{-\vert x-y\vert^2/4\kappa(t_2-s)}}
{(4\pi\kappa(t_2-s))^{n/2}}F(s,y)dyds=I_1(t_1,t_2,x)+I_2(t_1,t_2,x).
\enbal
\ee
The second term above satisfies the simple estimate
\begin{equation}\label{jul2618}
|I_2(t_1,t_2,x)|\le \|F\|_{L^\infty(\Rm\times\Rm^n)}|t_2-t_1|\le K|t_2-t_1|,
\end{equation}
thus we only need to look at $I_1$. 
For this term, we write, using the Newton-Leibniz formula in the $t$ variable, for a fixed $s\in[0,t_1]$
\be\label{jul2620}
\frac{e^{-\vert x-y\vert^2/4\kappa(t_2-s)}}{(4\pi\kappa(t_2-s))^{n/2}}
-\frac{e^{-\vert x-y\vert^2/4\kappa(t_1-s)}}{(4\pi\kappa(t_1-s))^{n/2}}=C_\kappa\int_{t_1}^{t_2}\frac{h(z)}
{(\tau-s)^{n/2+1}}d\tau,~~z=\frac{x-y}{\sqrt{\kappa(\tau-s)}},
\ee
with  an integrable function  
\be\label{oct2412}
h(z)=\Big(-\frac{n}2+\frac{\vert z\vert^2}{4}\Big)
e^{-\vert z\vert^2/4} .
\ee
Note that $h(z)$ has mean zero:
\be\label{oct2414}
\int_{\Rm^n}h(z)dz=0.
\ee
Thus, we have
\be
\bal
I_1(t_1,t_2,x)&= \int_0^{t_1}\int_{\Rm^n}\Big(
\frac{e^{-\vert x-y\vert^2/4\kappa (t_1-s)}}{(4\pi\kappa (t_1-s))^{n/2}}
-\frac{e^{-\vert x-y\vert^2/4\kappa(t_2-s)}}{(4\pi\kappa(t_2-s))^{n/2}}\Big)F(s,y)dyds\\
&=
C_\kappa\int_0^{t_1}\int_{\Rm^n}\int_{t_1}^{t_2}h\Big(\frac{x-y}{\sqrt{4\kappa(\tau-s)}}\Big)F(s,y) 
\frac{d\tau dyds}
{(\tau-s)^{n/2+1}} \\
&=C_\kappa\int_0^{t_1}\int_{\Rm^n}\int_{t_1}^{t_2}h(z)F(s,x-z\sqrt{4\kappa(\tau-s)}) 
\frac{d\tau dsdz}{\tau-s} .
\enbal
\ee
Since $h(z)$ has zero integral, this can be written as
\be\label{oct2416}
\bal
I_1(t_1,t_2,x)
&=C_\kappa\int_0^{t_1}\int_{\Rm^n}\int_{t_1}^{t_2}h(z)\big(F(s,x-z\sqrt{4\kappa(\tau-s)}) -F(s,x)\big)
\frac{d\tau dsdz}
{\tau-s} .
\enbal
\ee
We may now use (\ref{jul2212}) and (\ref{oct2414}) to split $I_1$ as 
\be\label{oct2417}
\bal
I_1(t_1,t_2,x)=I_0(t_1,t_2,x)+\bar I_1(t_1,t_2,x)-I_3(t_1,t_2,x),
\enbal
\ee
with
\be\label{oct2424}
\bal
I_0(t_1,t_2,x)
&=C_\kappa\int_0^{t_1/2}\int_{\Rm^n}\int_{t_1}^{t_2}h(z)
F(s,x-z\sqrt{4\kappa(\tau-s)})  
\frac{d\tau dsdz}
{\tau-s},
\enbal
\ee
and
\be\label{oct2421}
\bal
\bar I_1(t_1,t_2,x)
&=C_\kappa\int_{t_1/2}^{t_1}\int_{\Rm^n}\int_{t_1}^{t_2}h(z)
\big(\bar u(s,x-z\sqrt{4\kappa(\tau-s)}) -\bar u(s,x)\big)
\frac{d\tau dsdz}
{\tau-s},
\enbal
\ee
and
\be\label{oct2422}
\bal
I_3(t_1,t_2,x)
&=C_\kappa\int_{t_1/2}^{t_1}\int_{\Rm^n}\int_{t_1}^{t_2}h(z)
\big(u(s,x-z\sqrt{4\kappa(\tau-s)}) -u(s,x)\big)
\frac{d\tau dsdz}
{\tau-s} .
\enbal
\ee
The term $I_0$ is bounded in the same way as $I_2$:
\be\label{oct2425}
\bal
|I_0(t_1,t_2,x)|
&\le C\|F\|_{L^\infty}\int_0^{t_1/2}\int_{\Rm^n}\int_{t_1}^{t_2}|h(z)|
\frac{d\tau dsdz}{\tau-s}\le C\|F\|_{L^\infty}\int_0^{t_1/2}\int_{t_1}^{t_2} 
\frac{d\tau dsdz}{\tau-s}
\le C|t_2-t_1|,
\enbal
\ee
since $t_2\ge t_1\ge 2$.

As the function $u_0(x)$ is bounded, there exists
a constant $C$ so that
\be\label{oct814}
|\bar u(t,x)-\bar u(t,y)|\le \farc{C}{t^{1/2}}|x-y|,~~\hbox{for all $t\ge 1$ and $x,y\in\Rm^n$.} 
\ee
 We also have the uniform bound $0\le \bar u(t,x)\le K$ for all $t\ge 0$ and $x,y\in\Rm^n$.
This gives the following bound on $\bar I_1$
\be
\bal
&|\bar I_1(t_1,t_2,x)|
\le C\int_{t_1/2}^{t_1}\int_{\Rm^n}\int_{t_1}^{t_2}|h(z)|
\farc{C|z|(\tau-s)^{1/2}}{s^{1/2 } } 
\frac{d\tau dsdz}{\tau-s} 
=C\int_{t_1/2}^{t_1}\int_{t_1}^{t_2}\farc{d\tau ds}{s^{1/2}(\tau -s)^{1/2}}\\
&\le \farc{C}{t_1^{1/2}}\int_{t_1/2}^{t_1}\int_{t_1}^{t_2}\farc{d\tau ds}{(\tau -s)^{1/2}}\le
 \farc{C}{t_1^{1/2}}\int_{t_1/2}^{t_1}(\sqrt{t_2-s}-\sqrt{t_1-s})ds\\
 &\le \farc{C}{t_1^{1/2}}\int_{t_1/2}^{t_1}\farc{t_2-t_1}{\sqrt{t_2-s}+\sqrt{t_1-s}}ds
\le C|t_2-t_1|.
\enbal
\ee
Thus, our focus is now on the term $I_3$. We choose $q\in(0,1/10)$ and split this term as
\be
I_3(t_1,t_2,x)=C_\kappa(A_3(t_1,t_2,x)+R_3(t_1,t_2,x)),
\ee 
with
\be\label{oct2426}
\bal
A_3(t_1,t_2,x)
&= \int_{t_1/2}^{t_1}\int_{|z|\le (t_1-s)^q}\int_{t_1}^{t_2}h(z)
\big(u(s,x-z\sqrt{4\kappa(\tau-s)}) -u(s,x)\big)
\frac{d\tau dsdz}
{\tau-s}
\enbal
\ee
and
\be\label{oct2427}
\bal
R_3(t_1,t_2,x)
&= \int_{t_1/2}^{t_1}\int_{|z|\ge (t_1-s)^q}\int_{t_1}^{t_2}h(z)
\big(u(s,x-z\sqrt{4\kappa(\tau-s)}) -u(s,x)\big)
\frac{d\tau dsdz}
{\tau-s}.
\enbal
\ee
The term $R_3$ can be bounded as
\be
\bal
|R_3(t_1,t_2,x)|&\le C\|u\|_{L^\infty} \int_{t_1/2}^{t_1}
\int_{t_1}^{t_2}
e^{-c(t_1-s)^q}\frac{d\tau ds}{\tau-s}\\
&\le C\|u\|_{L^\infty} \int_{t_1/2}^{t_1}
e^{-c(t_1-s)^q}\big(\log(t_2-s)-\log(t_1-s)\big)ds\\
&\le C\|u\|_{L^\infty} \int_0^{\infty} e^{-cs^q}\big(\log(t_2-t_1+s)-\log s\big)ds.
\enbal
\ee
Thus, for any $m>0$ and $N>10$ we have
\be\label{oct2539}
\bal
&|R_3(t_1,t_2,x)|\le C\|u\|_{L^\infty} \int_0^{N} \big(\log(t_2-t_1+s)-\log s\big)ds+
C\|u\|_{L^\infty} \int_{N}^\infty e^{-cs^q}\big(\log(1+s)-\log s\big)ds\\
&\le \farc{C\|u\|_{L^\infty} }{N^m}
+ C\|u\|_{L^\infty}\int_0^N \big(\log(t_2-t_1+s)-\log s\big)ds=
\farc{C\|u\|_{L^\infty} }{N^m}\\
&+ {C\|u\|_{L^\infty} }
\Big((N+t_2-t_1)[\log(N+t_2-t_1)-1]-(t_2-t_1)(\log(t_2-t_1)-1)
-N(\log N-1)\Big)\\
&\le \farc{C\|u\|_{L^\infty} }{N^m}+ {C\|u\|_{L^\infty} }
\Big((N+t_2-t_1)\Big(\log N+C\farc{t_2-t_1}{N}-1\Big)-N(\log N-1)\\
&-(t_2-t_1)(\log(t_2-t_1)-1)\Big)\le  \farc{C\|u\|_{L^\infty} }{N^m}+ {C\|u\|_{L^\infty} }\Big((t_2-t_1)\log N+(t_2-t_1)|\log(t_2-t_1)|\Big).
 \enbal
\ee 
We used the fact that $0\le t_2-t_1\le 1$ above. 
Choosing $N=1/(t_2-t_1)$ gives
\be
\bal
|R_3(t_1,t_2,x)|&\le C_\alpha\|u\|_{L^\infty}|(t_2-t_1)^\alpha,
\enbal
\ee
for any $\alpha\in(0,1)$. 

The above analysis shows that 
the main task is to bound the term $A_3$ defined in (\ref{oct2426}). 
As a preliminary step,  note that for any fixed $T>0$  it suffices consider 
only~$t_1>10T$. This follows from the  local existence bounds already established in~\cite{CX1996} and~\cite{HLV}. 
In addition, observe that for any~$T>0$ fixed and $t_1>10T$ we have 
\be\label{jan1006}
\bal
&\tilde A_{3,T}(t_1,t_2,x)
:=\Big| \int_{t_1-T}^{t_1}\int_{|z|\le (t_1-s)^q}\int_{t_1}^{t_2}h(z)
\big(u(s,x-z\sqrt{4\kappa(\tau-s)}) -u(s,x)\big)
\frac{d\tau dsdz}
{\tau-s}\Big|\\
&\le C\int_{t_1-T}^{t_1}\int_{t_1}^{t_2}\farc{d\tau ds}{\tau -s}=C\int_{t_1-T}^{t_1}(\log(t_2-s)-\log(t_1-s))ds\\
&=C\int_0^T(\log(t_2-t_1+s)-\log s)ds\\
&=C\Big((t_2-t_1+T)(\log(t_2-t_1+T)-1)-(t_2-t_1)(\log(t_2-t_1)-1)-T(\log T-1))\\
&\le C\Big((t_2-t_1+T)(\log T+C\farc{t_2-t_1}{T})-1)-(t_2-t_1)(\log(t_2-t_1)-1)-T(\log T-1))\\
&\le C(1+\log T)|t_2-t_1|^\alpha,
\enbal
\ee
for any $\alpha\in(0,1)$. 

\subsection{A bound for the spatial variations of the fuel in the parabola}\label{sec:parabola}

The previous analysis shows that 
the main task is to bound the term $A_{3,T}$:
\be\label{oct2527}
\bal
A_{3,T}(t_1,t_2,x)
&= \int_{t_1/2}^{t_1-T}\int_{|z|\le (t_1-s)^q}\int_{t_1}^{t_2}h(z)
\big(u(s,x-z\sqrt{4\kappa(\tau-s)}) -u(s,x)\big)
\frac{d\tau dzds}
{\tau-s}.
\enbal
\ee
Recall that here we have $10T\le t_1\le t_2\le t_1+1$, and $T$ is fixed but sufficiently large, to be chosen. 
To estimate $A_{3,T}$, we need to show that the variations of the fuel concentration~$u(s,y)$ are small for most times~$s$ in the time interval $t_1/2\le s\le t_1$
inside the domain of integration in (\ref{oct2527}) that we will call a ``parabola",  centered at the point $x$ and the time $t_1$. 
We also define, for $A>10$ and for all $s\le t_1-T$ the region
\be\label{jan1216}
P^A_{s}=\big\{(s',y'):~0\le s'\le s,~|y'-x|\le A \sqrt{4\kappa}(t_1-s')^{{q}}{(t_1 + 1 - s')^{{1}/{2} } }\big\}
\ee
for the part of the parabola at times earlier than $s$ that we will refer to as sub-parabola. 
The constant~$A$ will be  
chosen to be large, and since $t_2 \le t_1 + 1$, it is clear that the points $(s,x)$ and~$(s, x-z \sqrt{4 \kappa (\tau - s)})$ in \eqref{oct2527} belong to $P^A_s$
if $|z|<(t_1-s)^q$.

An intuitive reason for the smallness of oscillations 
inside the parabola is that until the front arrives into this region, the fuel is
not consumed while its initial variations are quickly uniformized by the diffusion.  Once the front arrives, it spreads through the parabola with a linear speed,
consuming the fuel. The fuel variations during this period can not be controlled but this time interval is relatively short. After the front has passed,
there is no fuel left in the parabola, hence its variations are small again. Let us also comment that, according to Proposition~\ref{prop-oct2402}, once the
temperature is``large" at a given point, it never becomes ``too small" again. 

According to the above intuition, given a point $x\in\Rm^n$ and a time $t_1>10T$, with a sufficiently large but fixed $T>0$,  
we will consider  three cases: (1) the front ``has reached the parabola not long before $t_1$", (2) the front ``has reached the  
parabola long before $t_1$", and (3) the front ``has not reached the parabola by the time $t_1$".  

To quantify the above informal notions, let us first clarify the meaning of the phrase ``the front reached the parabola": 
if there is a time $s'\in [t_1/2,t_1-T]$ such that there exists a point~$x'$,
so that 
\be\label{oct2529}
|x-x'|\le {A} \sqrt{4\kappa}(t_1-s')^{{q}}{(t_1 + 1 - s')^{{1}/{2} }},
\ee
and the spatial average of the temperature around $x'$ is not small, in the sense that
\be\label{oct2528}
\int_{|y-x'|\le 1}v(s',y)dy\ge \farc{1}{(t_1-s')^m},
\ee
with $m>0$ to be chosen, we say that the $t_1$-front has reached the parabola centered at $x$ by the time~$s'$. 
Then, we define the ``front-reaching time" $s_0$ as the infimum of all times such that (\ref{oct2528}) holds:  
\be \label{feb0201}
s_0 := \inf \Big\{ s' \in [\frac{t_1}{2}, t_1 - T] ~\text{s.t. there is $x'$ so that (\ref{oct2529}) and (\ref{oct2528}) hold.}
\Big\}
\ee 
It follows that there is $x_0$ such that
\be\label{oct2529bis}
{|x-x_0|\le A \sqrt{4\kappa}(t_1-s_0)^{q} (t_1 + 1 - s_0)^{{1}/{2} },}
\ee
and
\be\label{oct2528bis}
\int_{|y-x_0|\le 1}v(s_0,y)dy\ge \farc{1}{(t_1-s_0)^m}.
\ee
Once the front-reaching time $s_0$ is set, we define the late-stage time $s_2$. 
Intuitively, it is the time when~$u$ is already decaying exponentially near the point $x$ due to the presence of front. First, we let~$s_1$ be
\be
s_1 := s_0 + m \log (t_1 - s_0) + C.
\ee
Corollary~\ref{cor-oct2402} implies that then the fuel concentration $u(s,y)$ obeys an upper bound
\be\label{oct2532}
u(s,y)\le Ke^{-\gamma_2(s-s_1)},~~\hbox{for all $s\ge s_1=s_0+m\log(t_1-s_0)+C$ 
and $|y-x_0|\le c_2(s-s_1)$.}
\ee
Here, we require $T>1$ so that Corollary \ref{cor-oct2402} applies.
Next, we define the late-stage time $s_2$ by
\be
s_2 := s_1 + \frac{1}{c_2} 10A \sqrt{4 \kappa} (t_1 - s_0)^q (t_1 - s_0 + 1)^{\frac{1}{2} },
\ee 
so that 
\be\label{oct2534}
c_2(s_2-s_1)= 10 A \sqrt{4\kappa}(t_1-s_0)^{q} (t_1 - s_0 + 1)^{\frac{1}{2} }.
\ee
holds. The intuition behind the choice of $s_2$ is clear: we know that $u$ is small 
near $x_0$ after the time~$s_1$. Then, by Corollary \ref{cor-oct2402}, after waiting for additional time $s_2 - s_1$, 
the smallness of $u$ is propagated and $u$ is exponentially decaying  throughout the parabola  after the time~$s_2$. 

Consider a time $s\in[s_2,t_1]$ and $(s,y)\in P^A_s$, that is, 
\be\label{oct2533}
|y-x|\le  A \sqrt{4\kappa}(t_1-s)^{q} (t_1 - s + 1)^{\frac{1}{2}}.
\ee
Then, due to (\ref{oct2529bis}), (\ref{oct2534}) and (\ref{oct2533}), we have
\be\label{oct2535}
|y-x_0|\le  A \sqrt{4\kappa}(t_1-s)^{q} (t_1 - s + 1)^{\frac{1}{2} } +A \sqrt{4\kappa}(t_1-s_0)^{q} (t_1 - s_0 + 1)^{\frac{1}{2} }
\le c_2(s_2-s_1)\le c_2(s-s_1).
\ee
We deduce from (\ref{oct2532}) and (\ref{oct2535}) that    the fuel concentration inside the parabola is, indeed, very small for $s>s_2$ by Corollary~\ref{cor-oct2402}:
\be\label{oct2536}
u(s,y)\le Ke^{-\gamma_2(s-s_1)},~~\hbox{for all $s\in [s_2,t_1]$ 
and $|y-x|\le A \sqrt{4\kappa}(t_1-s)^{q} (t_1 - s + 1)^{{1}/{2} } $.}
\ee
Let us comment that for convenience we require $T$ to be large enough, so that the function 
\[
t \rightarrow \frac{1}{2} t - m \log t - \frac{10 A \sqrt{4 \kappa} }{c_2} t^q (t+1)^{{1}/{2} } - C
\]
is positive for $t \ge T$. it follows, in particular, that 
\[
t_1 - s_2 > \frac{1}{2} (t_1 - s_0) \ge \frac{T}{2},
\]
that is, $s_2 < t_1 - {T}/{2}$. 

{\bf Brief outline of the strategy.}
With the above definitions of the times $s_0$ and $s_2$, we will split the integral \eqref{oct2527} into three pieces: 
first, the early times $s \in [t_1/2, s_0 - C \log (t_1 - s_0) - C]$. Here,~$C>0$ is a constant independent of $t_1, s_0$ that we will choose later.
The intermediate times cover the interval $s\in [s_0 - C\log (t_1 - s_0) - C, s_2]$, and late times include $s \in [s_2, t_1 - T]$. 
We will see that the integral over each time interval gives a contribution of the form $C|t_2 - t_1|$, where $C$ is independent of $t_1$, but for 
very different reasons.

First, for the early times $s$, by the backward in time estimate in Corollary \ref{cor-oct2404} and \eqref{feb0201}, 
we know that the 1-ball average of the fuel temperature $v$ near 
a point $x \in P^A_s$ decays exponentially. To get this decay,  
we needed to subtract from $s_0$ the additional term involving $\log(t_1 - s_0)$   in the definition of the upper bound for the early times. 
From Corollary \ref{cor-oct2506}, it follows that $v$ itself is small and moreover by the parabolic regularity theory, the gradient $\nabla v$ is 
also small in the smaller region~$P^{A'}_s$, with some $A'<A$. The application of the parabolic regularity theory explains why 
we needed to subtract the additional $O(1)$ term in the definition
of the upper bound for the early times.
Thus, during the early times and inside the slightly smaller parabola the fuel concentration 
$u$ is a solution to the heat equation with a small potential that also has a small gradient: see equation~\eqref{jan1202} below. 
Using Bernstein's method for the standard heat equation \cite{Krylov-book}, we then obtain a point-wise upper bound on
$\nabla u(s,x)$ inside an early times parabolic cylinder of size $R>0$. This bound is controlled by 
\[
\frac{1}{R^2} + R^2 \| \nabla v \|_{L^\infty_t L^2_y (Q_R^- (s,x) )}^2.
\]
By adjusting the size $R$ of the cylinder by the time $s$, to take advantage from smallness of $\nabla v$ in~$Q_R^- (s,x)$, we end up with 
the bound \eqref{jan1306} below for $\nabla u(s,y)$. This allows us to control the integral in (\ref{oct2527}) over the early times.

Second,  to control the integral in (\ref{oct2527}) over the late times, we use the exponential decay of $u$ in \eqref{oct2536}.

Finally,  on the intermediate time interval, neither $u$ nor $v$ are necessarily small. The key point is that if we only use the $L^\infty$ bound for the function
$u$,  the integral in \eqref{oct2527} is controlled by 
\[
\int  \big(\log (t_2 - s) - \log (t_1 - s)\big) ds \le \int \frac{t_2 - t_1}{t_1 - s} ds.
\]
However, in the intermediate time interval, the difference $t_1 -s$ varies within the range 
$$[(t_1 - s_0) - m \log (t_1 - s_0) - C - \frac{1}{c_2} 10 A \sqrt{4 \kappa} (t_1 - s_0)^q (t_1 - s_0 + 1)^{1/2}, (t_1 - s_0) + C \log (t_1 - s_0) + C].$$
 Since $T$ is sufficiently large and $t_1 - s_0 \ge T$, this means that integral of $(t_1 - s)^{-1}$ within this range of $s$ remains bounded uniformly in $t_1$. 
 
\commentout{
We need to distinguish two sub-cases: if $t_1/2<s_0\le t_1-T$, we will say that
the front has reached the parabola ``not long before $t_1$", and if $s_0=t_1/2$ we say that the front has reached the parabola ``long before~$t_1$".
If no time $s'$ as defined above exists in the time interval $[t_1/2,t_1-T]$,  we will say that the front has not reached the parabola by the time $t_1$. 
We will analyze these three cases separately.}

According to the above scenario, we 
estimate $A_{3,T}$ as 
\be\label{jan1327}
|A_{3,T}|\le A_e+A_m+A_\ell,
\ee
with
\be\label{oct2537}
\bal
A_e(t_1,t_2,x)
&= \int_{t_1/2}^{s_0 - C \log (t_1 - s_0) - C}\int_{|z|\le (t_1-s)^q}\int_{t_1}^{t_2}|h(z)|
\big|u(s,x-z\sqrt{4\kappa(\tau-s)}) -u(s,x)\big|
\frac{d\tau dzds}
{\tau-s},
\enbal
\ee
\be\label{oct2538}
\bal
A_m(t_1,t_2,x)
&= \int_{s_0 - C \log (t_1 - s_0) - C}^{s_2}\int_{|z|\le (t_1-s)^q}\int_{t_1}^{t_2}|h(z)|
\big|u(s,x-z\sqrt{4\kappa(\tau-s)}) -u(s,x)\big|
\frac{d\tau dzds}
{\tau-s},
\enbal
\ee
and
\be\label{oct2540}
\bal
A_\ell(t_1,t_2,x)
&= \int_{s_2}^{t_1-T/2}\int_{|z|\le (t_1-s)^q}\int_{t_1}^{t_2}|h(z)|
\big|u(s,x-z\sqrt{4\kappa(\tau-s)}) -u(s,x)\big|
\frac{d\tau dzds}
{\tau-s}.
\enbal
\ee
We have extended the domain of integration in $s$ in the term~$A_\ell$   
to the upper limit $s=t_1-T/2$ simply to avoid considering separately the cases $s_2<t_1-T$ and $s_2>t_1-T$,
while we do know that~$s_2<t_1-T/2$.

\commentout{
\subsubsection{The front has reached the parabola not long before the time $t_1$}
\label{sec:parab-reached}}
\commentout{
We first consider the most interesting case when $t_1/2<s_0\le t_1-T$. Informally, in that case we may think of $s_0$ as the first time reaction starts in the region around the point $x$. 
In that situation, we have equality in~(\ref{oct2528bis}):
\be\label{jan1002}
\int_{|y-x_0|\le 1}v(s_0,y)dy= \farc{1}{(t_1-s_0)^m}.
\ee
Corollary~\ref{cor-oct2402} implies that then the fuel concentration $u(s,y)$ obeys an upper bound
\be\label{oct2532}
u(s,y)\le Ke^{-\gamma_2(s-s_1)},~~\hbox{for all $s\ge s_1=s_0+m\log(t_1-s_0)+C$ 
and $|y-x_0|\le c_2(s-s_1)$.}
\ee
}
\commentout{Note that, as we have $s_0<t_1-T$, we know that $s_1<t_1-T/2$, as long as we choose $T$ sufficiently large, depending on the constants~$m$ and~$C>0$
in (\ref{oct2532}). }
 
\commentout{Next we take the  time $s_2\ge s_1$ such that 
\be\label{oct2534}
c_2(s_2-s_1)= 10\sqrt{4\kappa}(t_1-s_0)^{1/2+q}.
\ee
Once again, because $s_0<t_1-T$, we have $s_2<t_1-T/2$, as long as we choose $T$ sufficiently large, depending on $m$ and $C>0$. 
Physically, this is the time it takes the reaction to sweep across the region of integration in the definition of $A_3(t_1,t_2,x)$. 

Consider a time $s\in[s_2,t_1]$ and $y\in\Rm^n$ such that 
\be\label{oct2533}
|y-x|\le  \sqrt{4\kappa}(t_1-s)^{1/2+q}.
\ee
Then, due to (\ref{oct2529bis}), (\ref{oct2534}) and (\ref{oct2533}), we have
\be\label{oct2535}
|y-x_0|\le \sqrt{4\kappa}(t_1-s)^{1/2+q}+\sqrt{4\kappa}(t_1-s_0)^{1/2+q}
\le c_2(s_2-s_1)\le c_2(s-s_1).
\ee
We deduce from (\ref{oct2532}) and (\ref{oct2535}) that    the fuel concentration inside the parabola is, indeed, very small for $s>s_2$:
\be\label{oct2536}
u(s,y)\le Ke^{-\gamma_2(s-s_1)},~~\hbox{for all $s\in [s_2,t_1]$ 
and $|y-x|\le \sqrt{4\kappa}(t_1-s)^{1/2+q}$.}
\ee}
\commentout{
According to the above scenario, we 
estimate $A_{3,T}$ as 
\be\label{jan1327}
|A_{3,T}|\le A_e+A_m+A_\ell,
\ee
with
\be\label{oct2537}
\bal
A_e(t_1,t_2,x)
&= \int_{t_1/2}^{s_0}\int_{|z|\le (t_1-s)^q}\int_{t_1}^{t_2}|h(z)|
\big|u(s,x-z\sqrt{4\kappa(\tau-s)}) -u(s,x)\big|
\frac{d\tau dsdz}
{\tau-s},
\enbal
\ee
\be\label{oct2538}
\bal
A_m(t_1,t_2,x)
&= \int_{s_0}^{s_2}\int_{|z|\le (t_1-s)^q}\int_{t_1}^{t_2}|h(z)|
\big|u(s,x-z\sqrt{4\kappa(\tau-s)}) -u(s,x)\big|
\frac{d\tau dsdz}
{\tau-s},
\enbal
\ee
and
\be\label{oct2540}
\bal
A_\ell(t_1,t_2,x)
&= \int_{s_2}^{t_1-T/2}\int_{|z|\le (t_1-s)^q}\int_{t_1}^{t_2}|h(z)|
\big|u(s,x-z\sqrt{4\kappa(\tau-s)}) -u(s,x)\big|
\frac{d\tau dsdz}
{\tau-s}.
\enbal
\ee
We have extended the domain of integration in $s$ in the term~$A_\ell$   
to the upper limit $s=t_1-T/2$ simply to avoid considering separately the cases $s_2<t_1-T$ and $s_2>t_1-T$,
while we do know that~$s_2<t_1-T/2$.  
}

{\subsubsection{An estimate of the middle term $A_m$}}

For the ''middle" term $A_m$, we note that it follows from the definitions of the times $s_1$ and $s_2$ in~(\ref{oct2532}) and~(\ref{oct2534}) that 
for a sufficiently large $T>0$, we have
\be
\frac{t_1 - s_2}{(t_1 - s_0) + C \log (t_1 - s_0) + C } 
\in[1/3, 1]
\ee
since $t_1 - s_0 \ge T$.
\commentout{for any $\eps>0$,
if we choose $T>0$ sufficiently large, we have 
\be
\bal
1&\ge\farc{t_1-s_2}{t_1-s_0}=\farc{t_1-s_0-(s_2-s_1+s_1-s_0)}{t_1-s_0}\ge 1-\eps.
\enbal
\ee}
Therefore, we have the following estimate for $A_m$
\be\label{oct2602}
\bal
|A_m(t_1,t_2,x)|
&\le  \int_{s_0 {- C \log (t_1 - s_0) - C} }^{s_2} \int_{t_1}^{t_2}\frac{d\tau ds}
{\tau-s}=\int_{s_0 {- C \log (t_1 - s_0) - C} }^{s_2}[\log(t_2-s)-\log(t_1-s)]
\\
&\!\!\!\!\!=\int_{t_1-s_2}^{t_1-s_0 {+ C \log (t_1 - s_0) + C}}\big(\log(t_2-t_1+s)-\log s)ds
\le C\int_{t_1-s_2}^{t_1-s_0{+ C \log (t_1 - s_0) + C}}\farc{t_2-t_1}{s}ds\\
&\!\!\!\!\!=C(t_2-t_1)\log\Big(\frac{t_1-s_0{+ C \log (t_1 - s_0) + C}}{t_1-s_2}\Big)\le C(t_2-t_1).
\enbal
\ee

{\subsubsection{An estimate of the late term $A_\ell$}}

The ''late" term $A_\ell$ can be bounded using the upper bound (\ref{oct2536}):
\be\label{jan1328}
\bal
|A_\ell(t_1,t_2,x)|&\le C
\!\int_{s_2}^{t_1-T/2}\!\! \int_{t_1}^{t_2} 
e^{-\gamma_2(s-s_1)} 
\frac{d\tau dsdz}
{\tau-s}\le \! C\int_{s_{{2}}}^{t_1{-T/2}} e^{-\gamma_2(s-s_1)}\big(\log(t_2-s)-\log(t_1-s))ds\\
&=C\int_{s_{{2}}}^{t_1-T/2} e^{-\gamma_2(s-s_1)}\log\Big(\farc{t_2-s}{t_1-s}\Big)ds=
C\int_{s_{{2}}}^{t_1-T/2} e^{-\gamma_2(s-s_1)}\log\Big(1+\farc{t_2-t_1}{t_1-s}\Big)ds.
\enbal
\ee
As $0\le t_2-t_1\le 1$ and $t_1-s \ge T/2$, we can estimate the integral above as 
\be\label{jan1329}
\bal
|A_\ell(t_1,t_2,x)|&\le 
C\int_{s_{\red{2}}}^{t_1-T/2} e^{-\gamma_2(s-s_1)}\farc{t_2-t_1}{t_1-s}ds\le \farc{C}{T}|t_2-t_1|.
\enbal
\ee
\commentout{
We now split the last integral as in (\ref{oct2539}):
\be
\bal
|A_\ell(t_1,t_2,x)|\le C\int_{0}^{\infty} e^{-\gamma_2s}\big(\log(t_2-s-s_1)-\log(t_1-s-s_1))ds
\le C_\alpha(t_2-t_1)^\alpha\l
\enbal 
\ee
The last inequality above is obtained as in (\ref{oct2539}). }
\commentout{
For the "middle" term $A_m$, we note that it follows from  the definitions of the times $s_1$ and $s_2$ in~(\ref{oct2532}) and~(\ref{oct2534}) that for any $\eps>0$,
if we choose $T>0$ sufficiently large, we have 
\be
\bal
1&\ge\farc{t_1-s_2}{t_1-s_0}=\farc{t_1-s_0-(s_2-s_1+s_1-s_0)}{t_1-s_0}\ge 1-\eps.
\enbal
\ee
Therefore, we have
\be\label{oct2602}
\bal
|A_m(t_1,t_2,x)|
&\le  \int_{s_0}^{s_2} \int_{t_1}^{t_2}\frac{d\tau ds}
{\tau-s}=\int_{s_0}^{s_2}[\log(t_2-s)-\log(t_1-s)]
\\
&=\int_{t_1-s_2}^{t_1-s_0}\big(\log(t_2-t_1+s)-\log s)ds
\le C\int_{t_1-s_2}^{t_1-s_0}\farc{t_2-t_1}{s}ds\\
&=C(t_2-t_1)\log\Big(\frac{t_1-s_0}{t_1-s_2}\Big)\le C(t_2-t_1).
\enbal
\ee
}

\subsubsection{An estimate of the early term $A_e$}


The last and longer step is to estimate the early contribution $A_e$. To do this, we will estimate
the H\"older constant of the fuel concentration $u(t,x)$ inside the parabola for times 
$s$ in the interval
\[
t_1/2<s<s_0- C \log (t_1 - s_0) - C.
\]
The idea is as follows. We already know that the 
temperature $v(s,y)$ is small in all of the sub-parabola $P^{A}_{s_0}$. The parabolic regularity estimates
that we detail below will also imply that  the gradient $\nabla v(s,y)$ is small in the sub-parabola $P^{A'}_{s_0}$ for $A'<A$ sufficiently close to $A$. 
Therefore, the fuel concentration   is a solution to a parabolic equation of the form
\be
u_t=\Delta u-V(t,x)u,
\ee
with the potential $V(t,x)=g(v(t,x))$ that is very small and 
positive and also has a very small gradient in any sub-parabola~$P_s^{A'}$,
with $0\le s\le s_0$. 

In addition, if we take a time $s<s_0$, and a point $y\in P_{s_0}^{A'}$ that is not too close to the boundary of the parabola, then
we will be able to fit a parabolic cylinder~$Q_R(s,y)$ of a radius
\[
R\sim (t_1-s)^{1/2+q'},
\]
with some~$0<q'<q$, inside the parabola~$P_s^{A'}$. Hence, the fuel concentration~$u(t,x)$ will basically solve the heat equation inside 
a very large parabolic cylinder~$Q_R(s,y)$. This will show that~$|\nabla u(s,y)|\sim O(R^{-1})$,
and such bound will suffice for us to estimate $A_e$.  

We now follow the above outline. 
From the definition of $s_0$, we know that either $s_0 = {t_1}/{2}$, which is the case where $A_e = 0$, or $s_0 >  {t_1}/{2}$ and 
the local spatial averages of the temperature satisfy
\be
\int_{|y-z|\le 1}v(s_0,z)dz\le\farc{1}{|t_1-s_0|^m},~~\hbox{for each $y\in P^A_s$, that is, $|x-y|\le  A \sqrt{4\kappa}(t_1-s_0)^{q} (t_1 -s_0 + 1)^{\frac{1}{2} } $. } 
\ee
It follows from Corollary~\ref{cor-oct2404} that then  the spatial averages at the
earlier times satisfy
\be\label{oct2604}
\bal
\int_{|y-z|\le 1}v(s,z)dz &\le \frac{C}{|t_1-s_0|^m}  e^{-\gamma_3(s_0-s)},~~\hbox{for all $s\le s_0{- C \log (t_1 - s_0) - C}$ and } \\ 
&|x-y|\le {\farc{c_0}{2}(s_0-s)+ A \sqrt{4\kappa}(t_1-s_0)^{q} (t_1 - s_0 + 1)^{\frac{1}{2} } }.
\enbal
\ee
A key point is that Corollary~\ref{cor-oct2506} implies that then the point-wise
values of the temperature are also small, in the sense that they satisfy
\be\label{oct2702}
\bal
v(s,y) &\le  \frac{C}{|t_1-s_0|^{m(1-\eps)}}  e^{-\gamma_4(s_0-s)},~~\hbox{for all $s\le s_0 {- C \log (t_1 - s_0) - (C+4)}$ and } \\ 
&|x-y|\le \farc{c_0}{2}(s_0-s)+{A' \sqrt{4\kappa}(t_1-s_0)^{q} (t_1 - s_0 + 1)^{{1}/{2}}},
\enbal
\ee 
for some $A' \in (1, A)$. Note that we assume that $T$ is sufficiently large, so that 
\[
(A - A') \sqrt{4 \kappa} T^q (T+1)^{\frac{1}{2} } \ge C \ge 2,
\]
to guarantee that the points in $Q_2 (s,y)$ satisfy the conditions of \eqref{oct2604}.
Observe that, given~$c_0>0$, we can find $T>0$ so that if $t_1-s_0>T$ then the parabola $P^{A'}_s$ is contained in the domain of~\eqref{oct2702},
that is, we have
\be\label{jan1112}
{A' \sqrt{4\kappa}(t_1-s)^{q} (t_1 - s + 1)^{\frac{1}{2} } \le A' \sqrt{4\kappa}(t_1-s_0)^{q} (t_1 - s_0 + 1)^{\frac{1}{2} }+\farc{c_0}{2}(s_0-s),}
~~\hbox{for all $0\le s\le s_0$,}
\ee
since the function 
\[
r \rightarrow A' \sqrt{4 \kappa} r^q (r+1)^{{1}/{2} } - \frac{c_0}{2}r
\]
is decreasing for large enough $r$ as $q < \frac{1}{2} $.
It follows from (\ref{oct2702}) and 
(\ref{jan1112}) that the function $v(s,y)$ is very small in all of 
the sub-parabola ${P^{A'}_{s}}$ at the times $s$ earlier than $s_0 {- C \log (t_1 - s_0) - C}$, {after a suitable relabelling of the constant $C$}:
\be\label{jan1116}
\bal
v(s,y) &\le  \frac{C}{|t_1-s_0|^{m(1-\eps)}}  e^{-\gamma_4(s_0-s)},~~\hbox{for all $s\le s_0{- C\log (t_1 - s_0) - C}$ and } \\
& |x-y|\le {A' \sqrt{4\kappa}(t_1-s)^{q} (t_1 - s + 1)^{\frac{1}{2} }}.
\enbal
\ee 
Here, we again assume that $T$ is sufficiently large, so that the right-hand side of \eqref{jan1116} is less than~$1$. 

We will use the following   estimate for the solutions to the forced
heat equation (Theorem~12 in~\cite{Picard} and
Theorem 4.3 in~\cite{Kr-Saf}). Let $\phi(t,x)$ be a solution to 
\be\label{jan1120}
\phi_t=\Delta\phi+f(t,x).
\ee
Then, for any $t\ge 10$ and $x\in\Rm^n$  the function $\phi$ satisfies a H\"older estimate
\be\label{jan1121}
\|\phi\|_{C^{\alpha/2,\alpha}(Q_1\red{^-} (t,x))}\le C[\|\phi\|_{L^\infty(Q_2{^-}(t,x))}+\|f\|_{L^\infty(Q_2{^-}(t,x))}].
\ee
Therefore, the temperature $v(t,x)$ satisfies a bound
\be\label{jan1110}
\|v\|_{C^{\alpha/2,\alpha}(Q_1{^-})}\le 
C[\|v\|_{L^\infty(Q_2{^-})}+\|ug(v)\|_{L^\infty(Q_2{^-})}].
\ee
We deduce from (\ref{jan1116}), (\ref{jan1110}), and \eqref{oct2704} that $v(s,x)$ satisfies a uniform H\"older bound
\be\label{jan1118}
\bal
\|v\|_{C^{\alpha/2,\alpha}(Q_1{^-}(s,y))}&\le  \frac{C}{|t_1-s_0|^{m(1-\eps)}}  e^{-\gamma_4(s_0-s)},~~\hbox{for all $s\le s_0 {- C \log (t_1 - s_0) - C}$ and } \\
& |x-y|\le {A'' \sqrt{4\kappa}(t_1-s)^{q} (t_1 - s+1)^{\frac{1}{2} } }
\enbal
\ee 
{for $A'' \in (1, A')$ and sufficiently large $T$ so that for $s \le s_0 - C \log (t_1 - s_0) - C$ and $y \in P^{A''}_s$, we have
$\{z: |z-y| \le 2 \} \subset P^{A'}_s$.}
Similarly, it also follows from (\ref{jan1121})   that the fuel concentration $u(t,x)$ satisfies 
a uniform H\"older bound in the sub-parabola $P_{s_0}$:
\be\label{jan1123}
\bal
\|u\|_{C^{\alpha/2,\alpha}(Q_1^{-}(s,y))} &\le C,~~\hbox{for all $s\le s_0 {- C \log (t_1 - s_0) - C}$ and }\\
& |x-y|\le{A'' \sqrt{4\kappa}(t_1-s)^{q} (t_1 - s + 1)^{{1}/{2} }}.
\enbal
\ee 
We deduce from (\ref{jan1118}) and (\ref{jan1123}), as well as the assumption that $g$ is smooth and 
that $v(s,y)$ is small for $s \le s_0 - C \log (t_1 - s_0) - C$ and $y \in P^{A''}_s$, that the reaction rate also obeys a H\"older bound in $P_{s_0}$:
\be\label{jan1124}
\bal
\|ug(v)\|_{C^{\alpha/2,\alpha}(Q\red{^{-}}_1(s,y))}&\le  \frac{C}{|t_1-s_0|^{m(1-\eps)}}  e^{-\gamma_4(s_0-s)},~~\hbox{for all $s\le s_0{- C \log (t_1 - s_0) - C}$ and } \\
& |x-y|\le{A'' \sqrt{4\kappa}(t_1-s)^{q} (t_1 - s + 1)^{{1}/{2} }}.
\enbal
\ee 
Let us also recall Theorem~8.11.1 of~\cite{Krylov-book}. It says that a solution to (\ref{jan1120}) obeys a local estimate
\be\label{jan1122}
\|\phi\|_{C^{1+\alpha/2,2+\alpha}(Q^-_R)}\le C[\|f\|_{C^{\alpha/2,\alpha}(Q^-_{2R})}+\|\phi\|_{L^\infty(Q^-_{2R})}].
\ee 
This, together with (\ref{jan1116}) and (\ref{jan1124}) improves the bound (\ref{jan1118}) on the temperature to the following estimate in $P^{A''}_{s}$:
\be\label{jan1125}
\bal
\|v\|_{C^{1+\alpha/2,2+\alpha}(Q^-_{1/10}(s,y))} &\le  \frac{C}{|t_1-s_0|^{m(1-\eps)}}  e^{-\gamma_4(s_0-s)},~~\hbox{for all $s\le s_0{- C \log (t_1 - s_0) - C}$ and } \\
& |x-y|\le {A'' \sqrt{4\kappa}(t_1-s)^{q} (t_1 - s + 1)^{{1}/{2} } }.
\enbal
\ee 
We will mostly need the bound on $\nabla v(s,y)$ that is part of (\ref{jan1125}):
\be\label{jan1125bis}
\bal
\|\nabla v\|_{L^\infty(Q^-_{1/10}(s,y))}&\le  \frac{C}{|t_1-s_0|^{m(1-\eps)}}  e^{-\gamma_4(s_0-s)},~~\hbox{for all $s\le s_0{-C \log (t_1 - s_0) - C} $ and }\\
& |x-y|\le {A'' \sqrt{4\kappa}(t_1-s)^{q} (t_1 - s + 1)^{{1}/{2}} }.
\enbal
\ee

The next step is to translate the estimate (\ref{jan1125bis}) on the temperature 
into a small H\"older constant for the fuel concentration~$u(t,x)$ in $P_{s_0}$. Note that (\ref{jan1123}) gives a H\"older bound
on $u(t,x)$ that is not small and we need to improve that. To this end, let $\psi(t,x)$ be a positive solution to 
\be
\psi_t=\Delta \psi-V(t,x)\psi,
\ee
with a small non-negative potential $V(t,x)\ge 0$, in a parabolic cylinder 
\[
Q_R^- (s,x_1)=\{(t,x):~s-R^2< t <s,~~|x-x_1| <R\}.
\]  
We will assume without loss of generality that $s=0$
and $x_1=0$ and use the notation $Q_R=Q_R(0,0)$. 
Our goal is to get a bound on~$\nabla\psi(0,0)$ that would decay for $R$ large, similarly to a standard estimate 
for the heat equation, when $V\equiv 0$, as described, for instance,
in  Theorem 8.4.4 of~\cite{Krylov-book}. We will also assume that $\psi(t,x)$ is uniformly bounded:
\be\label{jan1220}
0\le\psi(t,x)\le K.
\ee
The constants below may depend on the constant $K$.

We follow the argument in~\cite{Krylov-book}. 
Let $\zeta_R(t,x)=\zeta(t/R^2,x/R)$, with $\zeta(t,x)$ a smooth cut-off function, such that~$\zeta(t,x)=1$ for~$-1/2\le t\le 0$
and $|x|\le 1/2$,  and $\zeta(s,y)=0$ if $s\le -2/3$ or~$|y|\ge 2/3$. We also define 
\be
w(t,x)=\zeta_R^2(t,x)|\nabla \psi(t,x)|^2+\mu|\psi(t,x)|^2,
\ee
with $\mu>0$ to be chosen. Then, we compute
\be
\bal
w_t&=\farc{2}{R^2}\zeta\zeta_t\Big(\farc{t}{R^2},\farc{x}{R}\Big)
|\nabla \psi(t,x)|^2+2\zeta_R^2 \nabla \psi\cdot\nabla \psi_t+2\mu \psi\psi_t\\
&=\farc{2}{R^2}\zeta\zeta_t\Big(\farc{t}{R^2},\farc{x}{R}\Big)|\nabla \psi(t,x)|^2+2\zeta_R^2 \nabla \psi\cdot\nabla(\Delta \psi)
-2\zeta_R^2\psi\nabla \psi\cdot\nabla V-2\zeta_R^2 V|\nabla\psi|^2\\
&+2\mu \psi\Delta \psi-2\mu V\psi^2,
\enbal
\ee
and (we use subscripts for partial derivatives in $x$ and sum over repeated indices) 
\be
\bal
\Delta w=\farc{2}{R^2}|\nabla \psi|^2(\zeta\Delta\zeta+|\nabla\zeta|^2)\Big(\farc{t}{R^2},\farc{x}{R}\Big)\!+2\zeta_R^2(\psi_{kj}\psi_{kj}
+\psi_k\Delta \psi_{k})\!+\!\farc{8}{R}
\zeta\zeta_k\psi_j\psi_{jk}+2\mu \psi\Delta \psi+2\mu |\nabla \psi|^2.
\enbal
\ee
This gives
\be\label{dec2702}
\bal
w_t-\Delta w&=|\nabla \psi|^2\Big(\farc{2}{R^2}\zeta\zeta_t-2\zeta_R^2 V-\farc{2}{R^2}(\zeta\Delta\zeta+|\nabla\zeta|^2)-2\mu\Big)
+2\zeta_R^2\big(\psi_k\psi_{kjj}-\psi_{kj}\psi_{kj}-\psi_k\psi_{kjj}\big)\\
&-2\zeta_R^2\psi \psi_kV_k-2\zeta_R^2 V|\nabla \psi|^2+2\mu \psi\Delta \psi-2\mu V\psi^2-\farc{8}{R}
\zeta\zeta_k\psi_j\psi_{jk}-2\mu \psi\Delta \psi \\
&\le |\nabla \psi|^2\Big(\farc{C}{R^2}-2\mu\Big)-2\zeta_R^2\psi \psi_kV_k-2\zeta_R^2\psi_{kj}\psi_{kj}-\farc{8}{R}
\zeta\zeta_k\psi_j\psi_{jk}-2\zeta_R^2 V|\nabla \psi|^2-2\mu V\psi^2\\
&= |\nabla \psi|^2\Big(\farc{C}{R^2}-2\mu\Big)-2\zeta_R^2\psi \psi_kV_k-2\Big(\zeta_R\psi_{kj}+\farc{2}{R}\zeta_k\psi_j\Big)^2
+\farc{8}{R^2}|\zeta_k|^2|\psi_j|^2  \\
& -2\zeta_R^2 V|\nabla \psi|^2-2\mu V\psi^2\\
&\le
 |\nabla \psi|^2\Big(\farc{C}{R^2}-2\mu\Big)-2\zeta_R^2\psi \psi_kV_k-2\zeta_R^2 V|\nabla \psi|^2-2\mu V\psi^2-2\Big(\zeta_R\psi_{kj}+\farc{2}{R}\zeta_k\psi_j\Big)^2.
\enbal
\ee
Let us denote
\be\label{dec2802}
g(t)=\|\nabla V(t,\cdot)\|_{L^\infty(B_R(0))}. 
\ee
As $V(t,x)\ge 0$ and $0\le \psi(t,x)\le C$, one can use (\ref{dec2702}) and (\ref{dec2802}) to write
\be\label{dec2706}
\bal
w_t-\Delta w 
&\le
 |\nabla \psi|^2\Big(\farc{C}{R^2}-2\mu\Big)+C|\nabla \psi||\nabla V|\le  |\nabla \psi|^2\Big(\farc{C}{R^2}-2\mu\Big)+CR^2|\nabla V|^2
 \le CR^2g^2(t),
\enbal
\ee
as long as $\mu>C/R^2$. Let us set 
\be
G(t)=CR^2\int_{-R^2}^tg^2(s)ds.
\ee
It follows that the function 
\be
\bar w(t,x)=w(t,x)-G(t)
\ee
satisfies the differential inequality 
\be
\bar w_t-\Delta\bar w\le 0,~~\hbox{ in $Q_R$.}
\ee
Therefore, $\bar w(t,x)$ attains its maximum on the parabolic boundary of $Q_R$. Moreover, $\zeta_R(t,x)=0$ if~$(t,x)$ lie on the parabolic boundary of $Q_R$,
and $\zeta_R(t,x)=1$ for $-R^2/2\le t\le 0$
and $|x|\le R/2$
Therefore, we have
\be
\bar w(t,x)\le \mu\|\psi\|_{L^\infty}^2\le \farc{C}{R^2}\|\psi\|_{L^\infty(Q_R)}^2,~~\hbox{ for all $(t,x)\in Q_{R/2}$}. 
\ee
We conclude that
\be\label{jan1210}
|\nabla \psi(0,0)|^2\le\farc{C\|\psi\|_{L^\infty}^2}{R^2}+G(0)=\farc{C\|\psi\|_{L^\infty}^2}{R^2}+CR^2\int_{-R^2}^0\|\nabla V(s,\cdot)\|_{L^\infty(Q_R)}^2ds.
\ee

We will use the above estimate as follows. Recall that $u(t,x)$ satisfies 
\be\label{jan1202}
u_t=\Delta u-g(v)u,
\ee
and that the nonlinearity $V(t,x)=g(v(t,x))$ satisfies the gradient bound (\red{\ref{jan1125bis}}) that implies, in particular, that
\be\label{jan1204}
\bal
|\nabla V(s,y)| &\le  \frac{C}{|t_1-s_0|^{m(1-\eps)}}  e^{-\gamma_4(s_0-s)},~~\hbox{for all {$(s,y) \in P^{A''}_{s_0 - C \log (t_1 - s_0) - C}.$}} 
\enbal
\ee 

Now, we focus  on the points 
\be \label{feb0501}
(s,y) \in P^1_{s_0 - C \log (t_1 - s_0) - C}, \frac{t_1}{2} \le s \le s_0 - C \log (t_1 - s_0) - C, 
\ee
so that  $|y-x| \le \sqrt{4\kappa} (t_1 - s)^q (t_1 + 1 - s)^{{1}/{2} }$ holds. These points are exactly those in the 
domain of integration for \eqref{oct2537}. We will choose the suitable size $R(s)$ of the parabolic cylinder so that~$Q_{R(s) }^- (s,y)$ is 
completely contained in $P^{A''}_{s_0 - C \log (t_1 - s_0) - C}$ and will deduce the smallness for $\nabla u$ from \eqref{jan1210}.

First, we fix $q'$ such that $0 < q' < q$ and choose $T>2$, so that $t_1 \ge 10 T > 20$. Next, for a time~$s \in \left [{t_1}/{2}, s_0 - C \log (t_1 - s_0) - C\right ]$, we define
\be \label{jan1221}
R(s) = \min ( \sqrt{4\kappa} (t_1 - s)^{q'} (t_1 + 1 -s)^{{1}/{2} }, (s-10)^{{1}/{2} } ).
\ee
Finally, we choose $A'' > 2$. Then for $(s,y)$ satisfying \eqref{feb0501}, we know that
\[
Q_{R(s)} ^- (s,y) \subset P^{A''}_{s_0 - C\log (t_1 -s_0) - C}.
\]
First, since $s \ge \frac{t_1}{2} >10$, $R(s) \ge 0$. Also, $s- s' < R^2 \le s-10$ implies that $s' > 10 > 0$.
Next, if $(s', y') \in Q_{R(s)} ^- (s,y)$, then 
\be
\bal
|y' - x| &\le |y' - y| + |y-x| \le R(s)+ \sqrt{4\kappa} (t_1 - s)^q (t_1 + 1 - s)^{\frac{1}{2} } \le A'' (t_1 - s)^q (t_1 + 1 - s)^{\frac{1}{2} }  \\
&\le A'' (t_1 - s')^q (t_1 + 1 - s')^{\frac{1}{2} }.
\enbal
\ee
Thus, $(s', y') \in P^{A''}_{s_0 - C \log (t_1 - s_0) - C}$ by the definition of the parabola \eqref{jan1216}.
 
\commentout{
Consider a time $s\in[t_1/2,s_0]$ and a point $y\in\Rm^n$ such 
that~$|y-x|\le (t_1-s)^{1/2+q'}$ with some~$0<q'<q$. 
Then, the parabolic cylinder 
\[
Q_R(s,y)=[s-R^2,s]\times B_R(y)
\]
is completely contained inside the sub-parabola $P_s$ as long as 
\be\label{jan1221}
0\le R\le R(s)=\min\big[(t_1-s)^{1/2+q'},\sqrt{s-10}\big].
\ee
That is, we have 
\be\label{jan1302}
\bal
&R(s)=(t_1-s)^{1/2+q'},~~\hbox{ for $s\ge t_c$},\\
&R(s)=\sqrt{s-10}, ~~\hbox{ for $t_1/2\le s\le t_c$},
\enbal
\ee
with the transition at the time 
\be\label{jan1320}
t_c\sim t_1^{1/(1+2q)}.
\ee
Note that it is possible that $t_c>s_0$ -- then we set $t_c=s_0$. }

We deduce then from (\ref{jan1210}) that  for all $(s,y)$ satisfying \eqref{feb0501}  we have
\be\label{jan1212}
\bal
|\nabla u(s,y)|&\le \farc{C}{R(s)}+CR(s)\Big(\int_0^s\frac{C}{|t_1-s_0|^{2m(1-\eps)}}  e^{-2\gamma_4(s_0-\tau)}d\tau\Big)^{1/2}\\
&\le
\farc{C}{R(s)}+\frac{C{(t_1-s)^{q'} (t_1 - s +1)^{{1}/{2} } }e^{-\gamma_4(s_0-s)}}{|t_1-s_0|^{m(1-\eps)}}.
\enbal
\ee
Note that if we choose $T$ sufficiently
large so that the function $t \rightarrow t^{q'+1} (t+1)^{\frac{1}{2} } e^{-\gamma_4 t }$ is decreasing for $t \ge T$,
then for all $s\le s_0\le t_1-T$ we have
\be
{(t_1-s)^{1+q'} (t_1 -s +1)^{{1}/{2} }e^{-\gamma_4(t_1-s)}\le (t_1-s_0)^{1+q'} (t_1 - s_0 + 1)^{{1}/{2} } e^{-\gamma_4(t_1-s_0)}.}
\ee
It follows that then
\be  {
\frac{(t_1-s)^{q'} (t_1 - s + 1)^{\frac{1}{2} }e^{-\gamma_4(s_0-s)}}{|t_1-s_0|^{m(1-\eps)}}\le 
\farc{(t_1-s_0)^{1+q'} (t_1 - s_0 + 1)^{\frac{1}{2} } }{(t_1-s_0)^{m(1-\eps)}(t_1-s)}\le\farc{C}{t_1-s},
}
\ee
as long as we choose $m>0$ sufficiently large. Here, the constant $C$ may depend on $T$. This gives the gradient 
bound
\be\label{jan1306}
\bal
|\nabla u(s,y)|&\le\farc{C}{R(s)}+\frac{C}{t_1-s}, \\&\hbox{ for all $t_1/2\le s \le s_0 {- C\log (t_1 - s) - C}$
and 
$|x-y|\le  (t_1-s)^{q}  {(t_1 - s + 1)^{{1}/{2} } }$.}
\enbal
\ee

We can now estimate $A_e$ as follows. Let us split this term as
\be\label{jan1308}
A_e=A_e^{(1)}+A_e^{(2)},
\ee
with 
\be\label{jan1310}
\bal
A_e^{(1)}(t_1,t_2,x)
&\!= \!\int_{t_1/2}^{s_0{- C \log (t_1 - s_0) - C}}
\!\!\int_{|z|\le (t_1-s)^{q-q'}}\!\int_{t_1}^{t_2}\!|h(z)|
\big|u(s,x-z\sqrt{4\kappa(\tau-s)}) -u(s,x)\big|
\frac{d\tau d{z}d{s}}
{\tau-s},
\enbal
\ee
and
\be\label{jan1312}
\bal
A_e^{(2)}(t_1,t_2,x)
&\!= \!\int_{t_1/2}^{s_0{- C \log (t_1 - s_0) - C}}
\!\!\!\int_{|z|\ge (t_1-s)^{q-q'}}\!\int_{t_1}^{t_2}|h(z)|
\big|u(s,x-z\sqrt{4\kappa(\tau-s)}) -u(s,x)\big|
\frac{d\tau d{z}d{s}}
{\tau-s}.
\enbal
\ee
The second term can be estimated as
\be\label{jan1314}
\bal
A_e^{(2)}(t_1,t_2,x)
&\le C\int_{t_1/2}^{s_0}\int_{t_1}^{t_2} e^{-  {c}(t_1-s)^{q-q'}}
\frac{d\tau ds}
{\tau-s}\le C|t_1-t_2|\int_{t_1/2}^{s_0} e^{-  {c}(t_1-s)^{q-q'}}
\farc{ds}{t_1-s}\\
&{ \le C |t_1 - t_2| \int_T ^\infty \frac{e^{-cr^{q-q'} }}{r} dr } \le C|t_1-t_2|.
\enbal
\ee
For the first term in (\ref{jan1308}) we use the gradient bound (\ref{jan1306}):
\be\label{jan1316}
\bal
A_e^{(1)}(t_1,t_2,x)
&\le C\int_{t_1/2}^{s_0}\int_{|z|\le (t_1-s)^{q-q'}}\int_{t_1}^{t_2}|h(z)||z|
\sqrt{\tau-s}\Big(\farc{1}{R(s)}+\frac{1}{t_1-s}\Big)
\frac{d\tau dsdz}
{\tau-s}\\
&\le C\int_{t_1/2}^{s_0}\int_{t_1}^{t_2}
\Big(\farc{C}{R(s)}+\frac{C}{t_1-s}\Big)
\frac{d\tau ds}
{\sqrt{\tau-s}}.
\enbal
\ee
We first note that
\be
\int_{t_1/2} ^{s_0} \int_{t_1} ^{t_2} \frac{1}{t_1 -s} \frac{d \tau d s}{\sqrt{\tau - s} } = 2 \int_{t_1/2} ^{s_0} \frac{(t_2 - t_1) }{(t_1 -s) (\sqrt{t_2 - s} + \sqrt{t_1 - s} )} ds \le C(t_2 - t_1)  \int_{T} ^\infty \frac{1}{r^{\frac{3}{2} } } dr \le C|t_1 - t_2|.
\ee
Also, since 
\be
\frac{1}{R(s) } \le \frac{1}{(s-10)^{\frac{1}{2} } } + \frac{1}{\sqrt{4\kappa} (t_1 - s)^{q'} (t_1 + 1 - s)^{\frac{1}{2} } } ,
\ee
we have
\be
\bal
\int_{t_1/2} ^{s_0} \int_{t_1} ^{t_2} \frac{d\tau ds} {R(s) \sqrt{\tau - s} } &\le \int_{t_1/2} ^{s_0} \int_{t_1} ^{t_2} \frac{d\tau ds} {(s-10)^{\frac{1}{2} } \sqrt{\tau - s} } + C \int_{t_1/2} ^{s_0} \int_{t_1} ^{t_2} \frac{d\tau ds} {  (t_1 - s)^{q'} (t_1 + 1 - s)^{\frac{1}{2} }  \sqrt{\tau - s} }  \\
&\le C |t_1 - t_2| \left ( \int_{t_1/2} ^{s_0} \frac{ds}{(s-10)^{\frac{1}{2} } (t_1 - s)^{\frac{1}{2} }} + \int_{t_1/2} ^{s_0} \frac{ds}{ (t_1 - s)^{q'+ \frac{1}{2} } (t_1 + 1 -s)^{\frac{1}{2} } }  \right ).
\enbal
\ee
By assuming $T>4$, we have that $s \ge t_1/2 \ge 5T > 20$, so that $s-10 > \frac{1}{2} s$. Then we have
\be
\int_{t_1/2} ^{s_0} \frac{ds}{(s-10)^{\frac{1}{2}} (t_1 - s)^{\frac{1}{2} } } \le  C \int_{t_1/2} ^{t_1} \frac{ds}{\sqrt{s (t_1 - s) } } = \int^1_{1/2} \frac{t_1 dr} { t_1 \sqrt{r(1-r) } } = C,
\ee
by change of variable $s = t_1 r$. On the other hand, we have
\be
\int_{t_1/2} ^{s_0} \frac{ds}{ (t_1 - s)^{q'+ \frac{1}{2} } (t_1 + 1 -s)^{\frac{1}{2} } } \le \int_{T} ^\infty \frac{dr}{r^{q'+1/2} (r+1)^{1/2} } \le C,
\ee
since $q'> 0$. Therefore, we have 
\be
A^{(1)}_e (t_1, t_2, x) \le C|t_1 - t_2|.
\ee

\commentout{
If $t_c$ given by (\ref{jan1320}) is larger than $s_0$ then the above expression is
\be\label{jan1321}
\bal
A_e^{(1)}(t_1,t_2,x)
&\le C\int_{t_1/2}^{s_0}\int_{t_1}^{t_2}
\Big(\farc{1}{\sqrt{s}}+\frac{1}{t_1-s}\Big)
\frac{d\tau ds}
{\sqrt{\tau-s}}\le C|t_1-t_2|\int_{t_1/2}^{s_0} 
\Big(\farc{1}{\sqrt{s}}+\frac{1}{t_1-s}\Big)
\frac{ds}{\sqrt{t_1-s}}\\
&\le C|t_1-t_2|\Big(1+\int_{t_1/2}^{t_1}\frac{ds}{\sqrt{s(t_1-s)}}\Big)
=C|t_1-t_2|\Big(1+\int_{1/2}^{1}\frac{ds}{\sqrt{s(1-s)}}\Big)=C|t_1-t_2|.
\enbal
\ee
We used above the fact that $s_0<t_1-T$. 
Finally, if $t_c<s_0$ 
then split this integral according to~(\ref{jan1302}):
\be\label{jan1318}
\bal
A_e^{(1)}(t_1,t_2,x)
&\le C\int_{t_1/2}^{t_c}\int_{t_1}^{t_2}
\Big(\farc{1}{\sqrt{s}}+\frac{1}{t_1-s}\Big)
\frac{d\tau ds}
{\sqrt{\tau-s}}+C\int_{t_c}^{s_0}\int_{t_1}^{t_2}
\Big(\farc{1}{(t_1-s)^{1/2+q'}}+\frac{1}{t_1-s}\Big)
\frac{d\tau ds}
{\sqrt{\tau-s}}\\
&\le C|t_1-t_2|\int_{t_1/2}^{t_c} 
\Big(\farc{1}{\sqrt{s}}+\frac{1}{t_1-s}\Big)
\frac{ds}
{\sqrt{t_1-s}}+C|t_1-t_2|
\int_{t_c}^{s_0} \farc{1}{(t_1-s)^{1/2+q'}} 
\frac{ds}
{\sqrt{t_1-s}}\\
&\le
C|t_1-t_2|.
\enbal
\ee
Here, the first integral was estimated exactly as in (\ref{jan1321}) and the second 
is finite because $s_0<t_1-T$. 

This shows that in the case when the time $s_0$ exists and is larger than $t_1/2$
we have the bound
\be\label{jan1322}
|A_3(t_1,t_2,x)|\le C|t_1-t_2|.
\ee

\subsubsection{The front has reached the parabola long before the time $t_1$
or has not reached it}

It remains to consider the cases when either the front has reached the parabola
long before $t_1$ -- recall that this means that $s_0=t_1/2$, or has not yet reached
it -- this means that for all $s\in[t_1/2,t_1]$ we have
\be\label{jan1325}
\int_{|y-x'|\le 1}v(s,y)dy\le \farc{1}{(t_1-s)^m}, 
~~~\hbox{ if $|x-x'|\le \sqrt{4\kappa}(t_1-s)^{1/2+q}$.}
\ee
 
Let us first assume that $s_0=t_1/2$. In that case, there exists $x_0$
such that 
\be
|x-x_0|\le\sqrt{4\kappa} (t_1-s_0)^{1/2+q},
\ee
and
\be\label{jan1326}
\int_{|y-x_0|\le 1}v(s,y)dy\ge \farc{1}{(t_1-s_0)^m}.
\ee 
Thus, the only terms in the decomposition (\ref{jan1327}) are $A_m$ and $A_\ell$ as the
early term $A_\ell$ is absent. The analysis of these terms in 
Section~\ref{sec:parab-reached}
(see (\ref{jan1328})-(\ref{jan1329}) for $A_\ell$ and (\ref{oct2602}) for $A_m$), 
only used the lower bound (\ref{jan1326}) and not the equality. Therefore, they
can proceed essentially verbatim the same, leading to (\ref{jan1329}) and
(\ref{oct2602}) also in this case.  

Finally, let us assume that the front has not yet arrived into the parabola,
so that (\ref{jan1325}) holds for all $t_1/2\le s\le t_1-T$. Then we can set
$s_0=t_1-T$ and only keep the early term $A_e$ in the decomposition (\ref{jan1327}).
The bounds for this term in Section~\ref{sec:parab-reached} 
only used the upper bound (\ref{jan1325}) at the time $s=s_0$. Therefore,
we may proceed exactly the same now as well. 

This finishes the proof of Lemma~\ref{lem-oct2402} and thus that 
of Theorem~\ref{thm:main} as well.
 
}


\begin{thebibliography}{99}
\bibliographystyle{plain}

%

\bibitem{Aronson-W}
D. Aronson and H. Weinberger, 
Multidimensional nonlinear diffusion arising in population genetics, Adv. Math. {\bf 30}, 1978, 33--76. 


%
 

  \bibitem{BHKR2005} H. Berestycki, F. Hamel, A. Kiselev and L. Ryzhik, Quenching and propagation in KPP reaction-diffusion equations with a heat loss, Arch. Ration. Mech. Anal. {\bf 178}, 2005, 57--80.
  
%
%

  
  

\bibitem{BN1992} J.Billingham and D.J.Needham,
The development of travelling waves in quadratic and cubic autocatalysis with unequal diffusion rates. III. Large time development in quadratic autocatalysis, Quart. Appl. Math., {\bf 50}, 1992, 343--372.

\bibitem{BD2006} M. Bisi and L. Desvillettes, From reactive Boltzmann equations to reaction-diffusion systems, J. Stat. Phys. {\bf 124}, 2006, 881-912.
\bibitem{Bramson1} 
M. D. Bramson, {Maximal displacement of branching Brownian motion}, 
{Comm. Pure Appl. Math.} {\bf 31}, 531--581, 1978.

\bibitem{Bramson2} 
M. D. Bramson, {\it Convergence of solutions of the Kolmogorov equation
 to travelling waves}, { Mem. Amer. Math. Soc.} {\bf 44}, 1983. 


%
%

  

%
 
 
%


%



\bibitem{CQ2008} X. Chen and Y. Qi, Propagation of local disturbances in reaction diffusion systems   modeling quadratic autocatalysis,
SIAM J. Appl. Math. {\bf 69}, 2008, 273--282. 

 
 
 \bibitem{CX1996} P. Collet and J. Xin,
Global existence and large time asymptotic bounds of $L^\infty$ solutions of thermal diffusive combustion systems on $\mathbb{R}^n$,  Ann. Sc. Norm. Sup. Pisa -- 
Classe di Scienze, Serie 4, {\bf 23}, 1996, 625--642. 

%
%
%
%


\bibitem{CS1977} D. E. Conway and J.A. Smoller, 
Diffusion and the predator-prey interaction. SIAM J. Appl. Math., {\bf 33}, 1977, 673--686.

%

\bibitem{DFT2017} L. Desvillettes, K. Fellner and B. Tang, Trend to equilibrium for reaction-diffusion systems arising from complex balanced chemical reaction networks, SIAM J. Math. Anal., {\bf 49}, 2017, 2666--2709.

 \bibitem{DFPV2007} L. Desvillettes, K. Fellner, M. Pierre, and J. Vovelle, 
About global existence for quadratic
systems of reaction--diffusion, Adv. Nonlinear Stud. {\bf{7}}, 2007, 491--511.

%
%
%
\bibitem{D2021}
A. Ducrot,  Spreading speed for a KPP type reaction-diffusion system with heat losses and fast decaying initial data,
Jour. Diff. Eqs {\bf 270}, 2021, 217--247.


 
%


%
\bibitem{FLT2018} K. Fellner, E. Latos and B. Tang,  Well-posedness and exponential equilibration of a volume-surface reaction-diffusion system with nonlinear boundary coupling. Ann. Inst. H. Poincar\'e   Anal. Non Lin\'eaire {\bf 35}, 2018,  643--673.



%
%
%
%
%
%

%
%
%
\bibitem{F1966} H. Fujita, On the blowing up of solutions of the Cauchy problem for $u_t=\Delta u+u^{1+\varepsilon}$,
J. Fac. Sci. Univ. Tokyo,  {\bf 13}, 1966, 109--124.
\bibitem{Graham}
C. Graham, Precise asymptotics for Fisher--KPP fronts, {Nonlinearity }
{\bf 32}, 2019, 1967--1998.  


 \bibitem{HR2005} F. Hamel and L. Ryzhik, Non-adiabatic KPP fronts with an arbitrary Lewis number, Nonlinearity {\bf 18}, 2005, 2881--2902.

\bibitem{HR2010} F. Hamel and  L. Ryzhik, Travelling waves for the thermodiffusive system with arbitrary Lewis numbers, Arch. Ration.
Mech. Anal. {\bf 195}, 2010, 923--952.

%
%


\bibitem{HLV} M. A. Herrero, A. A. Lacey, and J.L. Velazquez,
Global existence for reaction-diffusion systems modelling ignition, 
{Arch. Rational Mech. Anal} {\bf 142}, 1988, 219--251.

\bibitem{HMP1987} S.L. Hollis, R.H. Martin and M. Pierre, Global existence and boundedness in reaction--diffusion systems. SIAM J. Math. Anal. {\bf{18}}, 1987, 744--761.



%

%
%
%
%

 \bibitem{K2001} S. Kouachi, Existence of global solutions to reaction-diffusion systems via a Lyapunov
functional. Electron. J. Diff. Eqs, {\bf{68}}, 2001, 1--10.

%
%
%
%


\bibitem{Krylov-book}
N. Krylov, {\it Lectures on Elliptic and Parabolic Equations in H\"{o}lder Spaces}, Graduate Studies in Mathematics, 96. 
American Mathematical Society, Providence, RI, 2008. xviii+357 pp.

\bibitem{Kr-Saf}
N.V. Krylov and M.V. Safonov, A property of the solutions of parabolic equations 
with measurable coefficients, Izv. Akad. Nauk SSSR Ser. Mat., {\bf 44}, 1980,
161--175. 

\bibitem{Lieberman}
G.M. Lieberman, {\it Second Order Parabolic Differential Equations}, World Scientific,
1996.


\bibitem{MP1992}
R.H. Martin and M. Pierre,  
Nonlinear reaction-diffusion systems, in  
{\it Nonlinear equations in the applied sciences}, edited by W. F. Ames and C. Rogers, 363--398, 
Math. Sci. Engrg., 185, Academic Press, Boston, MA, 1992. 

%
\bibitem{MS1974} J. Maynard-Smith, {\it Models in Ecology}, Cambridge University Press, 1974.

\bibitem{MR2005} J.S. McGough and K. L. Riley, 
A priori bounds for reaction-diffusion systems arising
in chemical and biological dynamics, Appl. Math. Comput., {\bf 163}, 2005, 1--16.




%


\bibitem{OALD1997} H.G Othmer, F.R. Adler, M.A. Lewis and J. Dallon, eds, {\it Case Studies in Mathematical
Modeling-Ecology, Physiology and Cell Biology}, Prentice Hall, New Jersey, 1997.

\bibitem{Picard}
S. Picard, Notes on H\"older estimates for parabolic PDE, Preprint, 2019.

\bibitem{P2010} M. Pierre, Global existence in reaction-diffusion systems with control of mass: a survey, 
Milan J. Math. {\bf 78}, 2010, 417--455.

\bibitem{P1987} M. Pierre, An $L^1$-method to prove global existence in some reaction-diffusion systems. Contributions to nonlinear partial differential equations, Vol. II (Paris, 1985), 220--231, Pitman Res. Notes Math. Ser., {\bf 155}, Longman Sci. Tech., Harlow, 1987.
%
\bibitem{PS1997} M. Pierre and D. Schmitt, Blow up in reaction-diffusion systems with dissipation of
mass. SIAM J. Math. Anal. {\bf 28}, 1997, 259--269.


%
%
%

%
%


  
  

 \bibitem{S1993} E. M. Stein, {\it Harmonic Analysis: real-variable methods, orthogonality, and oscillatory integrals}, Princeton Mathematical Series, 43. Monographs in Harmonic Analysis, III. Princeton University Press, Princeton, NJ, 1993


  

  
  


%
%
 \end{thebibliography}
\end{document}